\documentclass{article}
\usepackage{amsmath}
\usepackage{amssymb}
\usepackage{graphicx}
\usepackage{epsfig}

\begin{document}

\renewcommand{\theequation}{\arabic{section}.\arabic{equation}}

\begin{center}
{\Large Sliding mode control of the Hodgkin--Huxley mathematical model}

\bigskip

Cecilia Cavaterra$^{1}$, Denis En\u{a}chescu$^{2}$, Gabriela Marinoschi$%
^{2,3}$

\bigskip

$^{1}$Dipartimento di Matematica ``F. Enriques'', Universit\`{a} degli Studi
di Milano

Via C. Saldini 50, 20133 Milano, Italy

and

Istituto di Matematica Applicata e Tecnologie Informatiche ``E. Magenes'',
CNR

Via Ferrata 1, 27100 Pavia, Italy

cecilia.cavaterra@unimi.it

\medskip

$^{2}$``Gheorghe Mihoc-Caius Iacob'' Institute of Mathematical Statistics
and Applied Mathematics of the Romanian Academy,

Calea 13 Septembrie 13, Bucharest, Romania

denaches@fmi.unibuc.ro

gabriela.marinoschi@acad.ro

\medskip

$^{3}$Research Group of the Project PN-III-P4-ID-PCE-2016-0372,

Simion Stoilow Institute of Mathematics

of the Romanian Academy, Bucharest, Romania

\bigskip
\end{center}

{\small \noindent Abstract. In this paper we deal with a feedback control
design for the action potential of a neuronal membrane in relation with the
non-linear dynamics of the Hodgkin-Huxley mathematical model. More exactly,
by using an external current as a control expressed by a relay graph in the
equation of the potential, we aim at forcing it to reach a certain manifold
in finite time and to slide on it after that. From the mathematical point of
view we solve a system involving a parabolic differential inclusion and
three nonlinear differential equations via an approximating technique and a
fixed point result. The existence of the sliding mode and the determination
of the time at which the potential reaches the prescribed manifold is proved
by a maximum principle argument. Numerical simulations are presented.}

\bigskip


MSC2010 Subject classification: 35K55, 35K57, 35Q92, 93B52, 92C30

\bigskip

\section{Introduction}

\setcounter{equation}{0}

The Hodgkin-Huxley (HH) model is the first complete mathematical model of
neuronal membrane dynamics explaining the ionic mechanisms determining the
initiation and propagation of action potentials in the squid giant axon. It
was successfully established\ in \cite{HH-1952} and since then it has become
a prototype model for all kinds of excitable cells, such as neurons and
cardiac myocytes. Detailed explanations of the biophysical process
illustrated by HH model can be found, e.g., in \cite{Best}, \cite{Chavarette}%
, besides the original work \cite{HH-1952}. In \cite{Chavarette} an analysis
of the non-linear dynamics in the Hodgkin-Huxley mathematical model showing
the existence of transient chaotic solutions in the model with their
original parameters, combined with the presentation of some modifications in
the dynamic system in order to become more realistic, has been done.

Many papers have been devoted to the mathematical analysis of this system
which exhibits a very complicated behavior. We confine ourselves to mention\
some fundamental mathematical works on the traditional Hodgkin-Huxley
equations: \cite{Cronin}, \cite{Evans-I}-\cite{Evans-Shenk}, \cite{Mascagni}%
, \cite{Fitzgibbon-93}. In the last one the existence of a unique classical
solution of the Hodgkin-Huxley system was proved. In the paper \cite{CC-MG},
the authors consider a singular perturbation of the Hodgkin-Huxley system
and study the associated dynamical system on a suitable bounded phase space,
when the perturbation parameter $\varepsilon $ (i.e., the axon specific
inductance) is sufficiently small, proving the existence of bounded
absorbing sets, of smooth attracting sets, as well as the existence of a
smooth global attractor.

From the mathematical point of view various properties of the dynamics of
the Hodgkin-Huxley vector field have been studied. Many studies in the
literature reveal bifurcations generated in the HH model, such as Hopf
bifurcation, period-double bifurcation and double cycle bifurcation (see
e.g., \cite{Che-2012} and the references there indicated). The HH model can
even exhibit a chaotic regime through a series of bifurcations. The
qualitative change of neuronal membrane potential from resting to repetitive
spiking, which is a characteristic behavior of this model, is of a
particular interest, because abnormal repetitive spiking are proper to
several neurological diseases. Consequently, much attention was directed to
provide mathematical results aiming to avoid instability around bifurcations
or to obtain desired dynamical behaviors which might be of help in the
development of the therapies of the diseases. For example, various dynamic
feedback control methods have been proposed to control the onset of Hopf
bifurcation in HH model, see, e.g., \cite{Che-2017} and the references there
indicated. We also cite the work \cite{Doruk-2010}, where the aim was to
develop a novel current control law with the purpose to stop the repetitive
firing caused by channel conductance deviations and the work \cite%
{Doruk-2013}, focusing on the simulation of the feedback controlled nerve
fiber stimulation where the behavior of the nerve fiber is manipulated by an
electrical field generator.

The Hodgkin-Huxley model introduced in \cite{HH-1952}, p. 522, eq. (29)
reads 
\begin{equation}
C_{M}\frac{dv}{dt}=\delta \frac{\partial ^{2}v}{\partial x^{2}}%
-g_{K}n^{4}(v-V_{K})-g_{Na}m^{3}h(v-V_{Na})-g_{l}(v-V_{l})+I_{C},\mbox{ in }%
Q,  \label{hh-1}
\end{equation}

\begin{equation}
\frac{dn}{dt}=\alpha _{n}(v)(1-n)-\beta _{n}(v)n,\mbox{ \ \ \ in }Q,
\label{hh-2}
\end{equation}%
\begin{equation}
\frac{dm}{dt}=\alpha _{m}(v)(1-m)-\beta _{m}(v)m,\mbox{ in }Q,  \label{hh-3}
\end{equation}%
\begin{equation}
\frac{dh}{dt}=\alpha _{h}(v)(1-h)-\beta _{h}(v)h,\mbox{ \ \ \ \ in }Q,
\label{hh-4}
\end{equation}%
where $v$ is the electrical potential in the nerve, $n,$ $m$, $h$ are the
proportions of the activating molecules of the potassium $(n),$ sodium
channels ($m)$ and of the inactivating molecules of the sodium channels $(h)$%
, respectively, $g_{K},$ $g_{Na},$ $g_{l}$ are the maximum conductances of
these ions, $V_{K},$ $V_{Na},$ $V_{l}$ are the constant equilibrium
potentials for these ions, $C_{M}$ is the membrane capacitance, $\delta =%
\frac{a}{2R_{2}}$ is a constant (depending on the fiber radius $a$ and the
specific resistance of the axoplasm $R_{2})$ and $I_{C}$ is the applied
current. Here $(t,x) \in Q:=(0,T)\times (0,L)$, where $x$ represents the
longitudinal distance along the axon and $t$ is time.

In (\ref{hh-2})-(\ref{hh-4}) $\alpha _{n},$ $\alpha _{m},$ $\alpha _{h},$ $%
\beta _{n},$ $\beta _{m},$ $\beta _{h}$ are nonlinear functions of $v $
defined as indicated, e.g., in \cite{Chavarette}, \cite{Doruk-2010}, \cite%
{Doruk-2013}, namely 
\begin{eqnarray}
\alpha _{n}(v) &=&\frac{0.1-0.01v}{\exp (1-0.1v)-1},\quad \mbox{ }\beta
_{n}(v)=0.125\exp \left( -\frac{v}{80}\right) ,  \label{alfa} \\
\alpha _{m}(v) &=&\frac{2.5-0.1v}{\exp (2.5-0.1v)-1},\quad \mbox{ }\beta
_{m}(v)=4\exp \left( -\frac{v}{18}\right) ,  \nonumber \\
\alpha _{h}(v) &=&0.07\exp \left( -\frac{v}{20}\right) , \quad \mbox{ }\beta
_{h}(v)=\frac{1}{\exp (3-0.1v)+1}.  \nonumber
\end{eqnarray}%
The values $\delta ,$ $C_{M},$ $g_{K},$ $g_{Na},$ $g_{l},$ $C_{M}$ are
positive numbers and $V_{K},$ $V_{Na},$ $V_{l}$ are real numbers.

This paper involves a new control approach, the sliding mode control, in
order to stabilize the membrane potential to a desired value. Sliding mode
control is\ an efficient tool for the stabilization of continuous or
discrete time systems. It consists in finding an appropriate control able to
constrain the evolution of the system in such a way to force it to reach a
manifold of a lower dimension, called the sliding manifold, in finite time,
and to keep it further sliding on this surface. Thus, our purpose is to
control the potential $v$ by means of a certain control $I_{C}$ in order to
force the potential to reach a prescribed value $v^{\ast }$ at a finite time 
$T^{\ast }$ and to keep this value for $t\geq T^{\ast }.$ The other state
variables $n,m,h$ will have after $T^{\ast }$ an evolution governed by their
equations in which $v$ takes the value $v^{\ast }.$ The principal advantage
of a sliding mode technique is that after some time the system evolves on a
manifold of lower dimension. For recent results regarding sliding mode
control for systems of parabolic equations we refer the reader to the papers 
\cite{BCGMR-sli}, \cite{Colli-Colturato}, \cite{CGMR-AMO}.

The objective is that $v$ reaches a constant value, in particular zero, and
to prove in this way the possibility to control the repetitive firing in
nerve fibers modeled by the Hodgkin-Huxley system. Even if a constant target
might be of main interest, the proof will be developed for a more general
case with $v^{\ast }$ dependending on time and space, which allows the
target to vary in time, being, for instance, periodic. To this end we
propose a relay feedback control of the form 
\begin{equation}
I_{C}(t,x)=-\rho \,\mbox{sign}\,(v(t,x)-v^{\ast }(t,x)),  \label{relay}
\end{equation}%
where the symbol sign denotes the multivalued function 
\begin{equation}
\mbox{sign\thinspace }r=\left\{ 
\begin{tabular}{ll}
$1,$ & $r>0$ \\ 
$\lbrack -1,1],$ & $r=0$ \\ 
$-1$, & $r<0$%
\end{tabular}%
\right.  \label{0-1}
\end{equation}%
and $\rho$ is a positive constant.

We rewrite (\ref{hh-1})-(\ref{hh-4}) in the following form 
\begin{equation}
C_{M}\frac{dv}{dt}=\delta \frac{\partial ^{2}v}{\partial x^{2}}%
-f_{1}(n,m,h)v+f_{2}(n,m,h)+I_{C},\mbox{ in }Q  \label{HH-new-1}
\end{equation}%
\begin{eqnarray}
\frac{dn}{dt} &=&-h_{1}^{n}(v)n+h_{2}^{n}(v),\mbox{ \ \ \ in }Q,  \nonumber
\\
\frac{dm}{dt} &=&-h_{1}^{m}(v)m+h_{2}^{m}(v),\mbox{ \ in }Q,
\label{HH-new-2} \\
\frac{dh}{dt} &=&-h_{1}^{h}(v)h+h_{2}^{h}(v),\mbox{ \ \ \ \ in }Q,  \nonumber
\end{eqnarray}%
where 
\begin{eqnarray}
f_{1}(n,m,h) &=&g_{K}n^{4}+g_{Na}m^{3}h+g_{l},  \label{fi-hi} \\
f_{2}(n,m,h) &=&g_{K}V_{K}n^{4}+g_{Na}V_{Na}m^{3}h+g_{l}V_{l},  \nonumber \\
h_{1}^{n}(v) &=&\alpha _{n}+\beta _{n},\quad \mbox{ }h_{1}^{m}(v)=\alpha
_{m}+\beta _{m},\quad \mbox{ }h_{1}^{h}(v)=\alpha _{h}+\beta _{h},  \nonumber
\\
h_{2}^{n}(v) &=&\alpha _{n},\quad \mbox{ }h_{2}^{m}(v)=\alpha _{m},\quad 
\mbox{ }h_{2}^{h}(v)=\alpha _{h}.  \nonumber
\end{eqnarray}%
The system is completed by homogeneous Neumann boundary conditions for $v,$%
\begin{equation}
\frac{\partial v}{\partial x}(t,0)=\frac{\partial v}{\partial x}(t,L)=0,%
\mbox{ for }t\in (0,T),  \label{bc}
\end{equation}%
since the membrane potential does not have a flux across the ends of the
fiber, and by initial conditions 
\begin{equation}
v(0,x)=v_{0},\mbox{ }n(0,x)=n_{0},\mbox{ }m(0,x)=m_{0},\mbox{ }h(0,x)=h_{0},
\quad x \in (0,L).  \label{ci}
\end{equation}

We shall approach this problem in two steps. First, as all equations for the
three components $n,$ $m$ and $h$ are similar, we shall consider a reduced
system formed only of two equations, one for the potential and the other for
only one ionic component, denoted generically by $w.$ This simplification
also occurs in the papers Fitzgibbon et al. (see \cite{Fitzgibbon-93}, \cite%
{Fitzgibbon-96}). In Section 2, we shall treat the simplified problem via an
approximating method, using a fixed point technique for proving the
existence of a solution to the system formed by the equation for the
membrane potential, with (\ref{relay}) replaced by involving the Yosida
approximation and one equation of the form (\ref{HH-new-2}). Suitable
estimates and compactness properties will allow to pass to the limit and to
prove an existence result for the non-approximated system in Theorem 2.1.
Then, the existence of the sliding mode will be provided in Theorem 2.2 by a
comparison argument. In Section 3, we shall extend the result to the
complete system (\ref{HH-new-1})-(\ref{HH-new-2}), by observing that it
follows as a consequence of the previous results for the simplified system.
The paper is concluded by numerical simulations intended to put into
evidence the sliding mode behavior of the solution.

\medskip

\noindent\textbf{Notation.} We denote%
\[
V=H^{1}(0,L)\subset H=L^{2}(0,L)\subset V^{\prime }=(H^{1}(0,L))^{\prime } 
\]%
where $V\subset H\subset V^{\prime }$ with compact injections. Moreover, we
define 
\[
W=\{y\in H^{2}(0,L);\mbox{ }y_{x}(0)=y_{x}(L)=0\}. 
\]

\begin{flushleft}
If $z\in L^{\infty }(X)$ the notation $\left\Vert z\right\Vert _{\infty }$
will stand for $\left\Vert z\right\Vert _{L^{\infty }(X)}$, where $X$ can be 
$\Omega ,$ or $Q.$ We denote by $C,$ $C_{i},$ $i=1,2,...$ some constants
depending on problem parameters, sometimes explicitly indicated in the
argument. For the sake of simplicity we shall write $v_{t},$ $v_{x},$ $%
v_{xx} $ instead of $\frac{d v}{d t},$ $\frac{\partial v}{\partial x},$ $%
\frac{\partial ^{2}v}{\partial x^{2}}$ and similarly for the other functions.
\end{flushleft}

\section{The simplified system}

\setcounter{equation}{0}

Let us consider the system for the potential $v$ and the concentration $w,$
coupled with a set of homogeneous Neumann boundary conditions for the
potential and of initial data%
\begin{equation}
v_{t}-\delta v_{xx}+f_{1}(w)v+\rho \mbox{\thinspace sign\thinspace }%
(v-v^{\ast })\ni f_{2}(w),\mbox{ in }Q,  \label{1}
\end{equation}%
\begin{equation}
w_{t}=-h_{1}(v)w+h_{2}(v),\mbox{ \ \ \ \ \ in }Q,  \label{2}
\end{equation}%
\begin{equation}
v_{x}(t,0)=v_{x}(t,L)=0,\mbox{ \ \ \ \ \ \ \ \ \ \ \ in }(0,T),  \label{3}
\end{equation}%
\begin{equation}
v(0,x)=v_{0},\mbox{ }w(0,x)=w_{0},\mbox{ \ \ \ \ \ in }(0,L).  \label{4}
\end{equation}%
The desired final value to be obtained is the time and space dependent
function $v^{\ast }$. Here the value $C_{M}$ is considered for simplicity
equal to 1.

Taking into account the general considerations presented in the introduction
on the expressions of the functions occurring in the Hodgkin-Huxley model,
we shall assume the following properties:

\begin{itemize}
\item[(i)] the functions $f_{i}$ and $h_{i}$, for $i=1,2,$ are locally
Lipschitz continuous, that is, for any $M$ positive, and for any $r$, $%
\overline{r}\in \mathbb{R}$, $\left\vert r\right\vert \leq M,$ $\left\vert 
\overline{r}\right\vert \leq M,$ there exist $L_{f_{i}}(M)$ and $%
L_{h_{i}}(M) $ positive, such that 
\begin{equation}
\left\vert f_{i}(r)-f_{i}(\overline{r})\right\vert \leq
L_{f_{i}}(M)\left\vert r-\overline{r}\right\vert ,\mbox{ \ }\left\vert
h_{i}(r)-h_{i}(\overline{r})\right\vert \leq L_{h_{i}}(M)\left\vert r-%
\overline{r}\right\vert;  \label{5}
\end{equation}

\item[(ii)] there exists $a>0$ such that 
\begin{equation}
0<a\leq f_{1}(r),\mbox{ }0<h_{i}(r),\mbox{ }i=1,2,\mbox{ for all }r\in 
\mathbb{R};  \label{5-0}
\end{equation}

\item[(iii)] 
\begin{equation}
v_{0}\in L^{2}(0,L),\mbox{ }w_{0}\in L^{2}(0,L);  \label{5-1}
\end{equation}

\item[(iv)] 
\begin{equation}
v^{\ast }\in W^{1,\infty}(0,T;H) \cap L^{\infty}(0,T;W).  \label{6}
\end{equation}
\end{itemize}

\medskip

\noindent \textbf{Definition 2.1.} We call a solution to system (\ref{1})-(%
\ref{4}) a pair $(v,w)$%
\[
v\in C([0,T];H)\cap L^{2}(0,T;V)\cap W^{1,2}(0,T;H)\cap L^{\infty }(Q),\mbox{%
\ }w\in L^{\infty }(Q)\cap W^{1,2}(0,T;H) 
\]
which satisfies 
\begin{eqnarray}
&&\int_{0}^{T}\int_{0}^{L}(v_{t}\psi +\delta \nabla v\cdot \nabla \psi
+f_{1}(w)v\psi +\rho \zeta \psi )dxdt  \label{7} \\
&=&\int_{0}^{T}\int_{0}^{L}f_{2}(w)\psi dxdt,\mbox{ for all }\psi \in
L^{2}(0,T;V),\mbox{ }  \nonumber \\
\zeta (t,x) &\in &\mbox{sign}\,(v(t,x)-v^{\ast }(t,x))\mbox{ a.e. }(t,x)\in
(0,T)\times (0,L),  \nonumber
\end{eqnarray}%
\begin{equation}
\int_{0}^{T}\int_{0}^{L}(w_{t}\phi +h_{1}(v)w\phi
)dxdt=\int_{0}^{T}\int_{0}^{L}h_{2}(v)\phi dxdt,\mbox{ for all }\phi \in
L^{2}(0,T;H),  \label{7-0}
\end{equation}%
together with the initial conditions (\ref{4}).

\medskip

We observe that by hypotheses (i) and (\ref{6}) it follows that $h_{i}(v)$
and $f_{1}(w)$ belong to $L^{\infty }(Q)$ and so the integrals containing
these functions make sense.

\medskip

\noindent \textbf{Theorem 2.1.} \textit{Let} (i)-(iv)\textit{\ hold. Assume
that} 
\begin{equation}
\sup_{r\in \mathbb{R}}\frac{h_{2}(r)}{h_{1}(r)}=: w_{M} \in \mathbb{R},\mbox{
}  \label{9-2}
\end{equation}%
\textit{and consider}%
\begin{equation}
v_{0}\in V,  \label{9}
\end{equation}%
\begin{equation}
w_{0}\in C([0,L]),\mbox{ \ }0\leq w_{0}(x)\leq w_{M}\mbox{ \textit{for all} }%
x\in \lbrack 0,L],  \label{9-1}
\end{equation}%
\begin{equation}
v^{\ast }\in W^{1,\infty }(0,T;H)\cap L^{\infty }(0,T;W).  \label{9-3}
\end{equation}%
\textit{Then, problem} (\ref{1})-(\ref{4}) \textit{has a unique solution,
with the further regularity} 
\begin{equation}
v\in L^{\infty }(0,T;V)\cap L^{2}(0,T;W),\mbox{ }w\in C([0,T];C[0,L])\cap
W^{1,\infty }(0,T;H).  \label{10}
\end{equation}

\medskip

\noindent \textbf{Proof. }We shall consider a regularized problem and prove
that it has a unique solution by applying the Banach fixed point theorem.
Then, we shall pass to the limit to recover the solution to (\ref{1})-(\ref%
{4}).

Let $\varepsilon $ be positive and introduce the Yosida approximation of the
sign operator, 
\[
\mbox{sign}_{\varepsilon }v=\frac{1}{\varepsilon }(I-(I+\varepsilon \mbox{%
sign})^{-1})v, 
\]%
and the approximating system 
\begin{equation}
(v_{\varepsilon})_{t}-\delta
(v_{\varepsilon})_{xx}+f_{1}(w_{\varepsilon})v_{\varepsilon}+\rho \,\mbox{%
sign}_{\varepsilon }(v_{\varepsilon}-v^{\ast })=f_{2}(w_{\varepsilon}),\mbox{
in }Q,  \label{1-1}
\end{equation}%
\begin{equation}
(w_{\varepsilon})_{t}=-h_{1}(v_{\varepsilon})w_{\varepsilon}+h_{2}(v_{%
\varepsilon}),\mbox{ \ \ \ \ \ \ \ \ \ \ \ \ \ \ \ \ \ \ \ \ \ \ \ \ \ \ \ \
\ in }Q,  \label{1-1bis}
\end{equation}%
\begin{equation}
(v_{\varepsilon})_{x}(t,0)=(v_{\varepsilon})_{x}(t,L)=0,\mbox{ \ }\ \ \ \ \
\ \ \ \ \ \ \ \ \ \ \ \ \ \ \ \ \ \ \ \ \ \ \ \ \ \ \ \ \ \mbox{in }(0,T)%
\mbox{,}  \label{3-1}
\end{equation}%
\begin{equation}
v_{\varepsilon}(0,x)=v_{0},\mbox{ }w_{\varepsilon}(0,x)=w_{0},\mbox{ \ \ \ \
\ \ \ \ \ \ \ \ \ \ \ \ \ \ \ \ \ \ \ \ \ \ \ \ \ in }(0,L).  \label{4-1}
\end{equation}%
Let $R$ be a positive\ value, which will be later specified, and let us
introduce the set 
\begin{eqnarray*}
\mathcal{M} &=&\left\{ (v,w)\in C([0,T];H)\times C([0,T];H);\right. \mbox{ }%
v\in L^{\infty }(0,T;V),\mbox{ }w\in L^{\infty }(Q), \\
&&\mbox{ \ \ \ \ \ \ \ \ \ \ \ \ \ \ }\left. \left\Vert v\right\Vert
_{L^{\infty }(0,T;V)}\leq R,\mbox{ }\left\Vert w\right\Vert _{L^{\infty
}(Q)}\leq w_{M}\right. \},
\end{eqnarray*}%
which obviously is a closed subset of $C([0,T];H)\times C([0,T];H)$. Also, ${%
\mathcal{M}}$ is a metric space with the metric $d_{\mathcal{M}} ((u,w),(%
\overline{u},\overline{w}))=\left\Vert u-\overline{u}\right\Vert
_{C([0,T];H)}+\left\Vert w-\overline{w}\right\Vert _{C([0,T];H)}.$ We shall
apply the Banach fixed point theorem in $\mathcal{M}.$

We fix $(\overline{v},\overline{w})\in \mathcal{M}$ and consider the system%
\begin{equation}
(y_{\varepsilon})_{t}-\delta (y_{\varepsilon})_{xx}+f_{1}(\overline{w}%
)y_{\varepsilon}=f_{2}(\overline{w})-\rho \,\mbox{sign}_{\varepsilon }(%
\overline{v}-v^{\ast }),\mbox{ in }Q,  \label{11}
\end{equation}%
\begin{equation}
(z_{\varepsilon})_{t}=-h_{1}(\overline{v})z_{\varepsilon} + h_{2}(\overline{v%
}),\mbox{ \ \ \ \ \ \ \ \ \ \ \ \ \ \ \ \ \ \ \ \ \ \ \ \ \ \ \ \ \ \ in }Q,
\label{12}
\end{equation}%
\begin{equation}
(y_{\varepsilon})_{x}(t,0)=(y_{\varepsilon})_{x}(t,L)=0,\mbox{ \ }\ \ \ \ \
\ \ \ \ \ \ \ \ \ \ \ \ \ \ \ \ \ \ \ \ \ \ \ \ \ \ \ \ \mbox{in }(0,T),
\label{13}
\end{equation}%
\begin{equation}
y_{\varepsilon}(0,x)=v_{0},\mbox{ }z_{\varepsilon}(0,x)=w_{0},\mbox{ \ \ \ \
\ \ \ \ \ \ \ \ \ \ \ \ \ \ \ \ \ \ \ \ \ \ \ \ \ in }(0,L).  \label{14}
\end{equation}%
Since%
\[
\left\vert f_{i}(\overline{w})\right\vert \leq \left\vert
f_{i}(0)\right\vert +L_{f_{i}}(w_{M})\left\vert \overline{w}\right\vert ,%
\mbox{ }\left\vert h_{i}(\overline{v})\right\vert \leq \left\vert
h_{i}(0)\right\vert +L_{h_{i}}(R)\left\vert \overline{v}\right\vert ,\mbox{ }%
i=1,2, 
\]%
with $\overline{v}\in L^{\infty }(0,T;C[0,L])$ (indeed $V\subset C[0,L]$ in
a one-dimensional space), we get $f_{i}(\overline{w})\in L^{\infty }(Q)$ and 
$h_{i}(\overline{v})\in L^{\infty }(Q).$ Let us set 
\begin{equation}
f_{iM}:=\left\vert f_{i}(0)\right\vert +L_{f_{i}}(w_{M})w_{M},\mbox{ }%
h_{iR}:=\left\vert h_{i}(0)\right\vert +L_{h_{i}}(R)R,\mbox{ }i=1,2.
\label{11-1}
\end{equation}%
We note that $f_{iM}$ depend on $w_{M}$ while $h_{iR}$ depend on $R,$ $i=1,2$
and take 
\begin{equation}
R\geq \overline{C},  \label{11-0}
\end{equation}%
where $\overline{C}$ is a constant depending on the problem parameters and
the initial datum for $v_{\varepsilon}$ and will be given below.

Next, we define $\Psi :\mathcal{M}\rightarrow L^{2}(0,T;H)\times
L^{2}(0,T;H),$ by $\Psi (\overline{v},\overline{w})=(y_{\varepsilon
},z_{\varepsilon })$ the solution to (\ref{11})-(\ref{14}) and prove further
that $\Psi (\mathcal{M})\subset \mathcal{M}$ and that $\Psi $ is a
contraction.

By (\ref{12}) we have 
\begin{equation}
z_{\varepsilon }(t,x)=e^{-\int_{0}^{t}h_{1}(\overline{v}(\sigma ,x))d\sigma
}w_{0}(x)+\int_{0}^{t}e^{-\int_{s}^{t}h_{1}(\overline{v}(\sigma ,x))d\sigma
}h_{2}(\overline{v}(s,x))ds.  \label{14-1}
\end{equation}%
It is immediately seen that $z_{\varepsilon }\in C([0,T];C[0,L])$. Indeed $%
(t,x)\rightarrow \int_{s}^{t}h_{1}(\overline{v}(\sigma ,x))d\sigma $ and $%
(t,x)\rightarrow \int_{0}^{t}e^{-\int_s^t h_1(\overline v(\sigma,
x))d\sigma} h_{2}(\overline{v}(\sigma ,x))d\sigma $, for $0\leq s<t\leq T $,
are continuous on $[0,T]$ and $w_{0}\in C([0,L])$. By (\ref{5-0}), (\ref{9-1}%
), (\ref{9-2}) it follows that $z_{\varepsilon }(t,x)\geq 0$ and 
\begin{eqnarray*}
\left\vert z_{\varepsilon }(t,x)\right\vert &\leq
&w_{M}+w_{M}\int_{0}^{t}e^{-\int_{s}^{t}h_{1}(\overline{v}(\sigma
,x))d\sigma }h_{1}(\overline{v}(s,x))ds \\
&=&w_{M}-w_{M}\left( 1-e^{-\int_{0}^{t}h_{1}(\overline{v}(\sigma ,x))d\sigma
}\right) \leq w_{M},\mbox{ for all }t\in \lbrack 0,T]\times \lbrack 0,L],
\end{eqnarray*}%
hence 
\begin{equation}
\left\Vert z_{\varepsilon }\right\Vert _{C([0,T];C[0,L])}\leq w_{M}.
\label{15}
\end{equation}%
Moreover, by (\ref{12}) we see that 
\begin{equation}
\left\Vert (z_{\varepsilon })_{t}\right\Vert _{C([0,T];C[0,L])}\leq
h_{1R}w_{M}+h_{2R}.  \label{15-0}
\end{equation}

In order to deal with the parabolic problem (\ref{11}), (\ref{13}), (\ref{14}%
) we introduce the linear time dependent operator $A(t):V\rightarrow
V^{\prime },$ 
\[
\left\langle A(t)y,\psi \right\rangle _{V^{\prime },V}=\int_{0}^{L}(\delta
y_{x}\psi _{x}+f_{1}(\overline{w}(t,x))y\psi )dx,\mbox{ for all }\psi \in V 
\]%
and write the equivalent Cauchy problem%
\begin{eqnarray}
\frac{dy_{\varepsilon }}{dt}(t)+A(t)y_{\varepsilon }(t) &=&f_{2}(\overline{w}%
(t))-\rho \,\mbox{sign}_{\varepsilon }(\overline{v}(t)-v^{\ast }(t)),\mbox{
a.e. }t\in (0,T),  \label{17} \\
y_{\varepsilon }(0) &=&v_{0}.  \nonumber
\end{eqnarray}%
The operator $A(t)$ has the properties 
\[
\left\Vert A(t)y\right\Vert _{V^{\prime }}\leq \max \{f_{1M},\delta
\}\left\Vert y\right\Vert _{V},\mbox{ }\left\langle A(t)y,y \right\rangle
_{V^{\prime },V}\geq \min \{a,\delta \}\left\Vert y\right\Vert _{V}^{2}, 
\]%
and so by the Lions theorem (see \cite{Lions-1969}, p. 162), the Cauchy
problem has a unique solution $y_{\varepsilon }\in W^{1,2}(0,T;V^{\prime
})\cap L^{2}(0,T;V)\cap C([0,T];H).$ The solution satisfies a first
estimate, obtained by testing (\ref{17}) by $y_{\varepsilon }(t)$ in $H$ and
then integrating over $(0,t)$%
\begin{equation}
\left\Vert y_{\varepsilon }(t)\right\Vert _{H}^{2}+\int_{0}^{t}\left\Vert
y_{\varepsilon }(s)\right\Vert _{V}^{2}ds\leq \frac{1}{\delta _{1}}\left(
\left\Vert v_{0}\right\Vert _{H}^{2} + 2(f_{2M}^{2}+\rho ^{2})LT\right) e^{%
\frac{2}{\delta_1}t},\mbox{ for all }t\in \lbrack 0,T],  \label{19}
\end{equation}%
where\ $\delta _{1}=\min \{1,2a,2\delta \}.$

We calculate a second estimate, by multiplying formally (\ref{17}) in $H$ by 
$-(y_{\varepsilon })_{xx}(t)$ and then integrating over $(0,t).$ We get%
\begin{eqnarray*}
&&\frac{1}{2}\left\Vert \nabla y_{\varepsilon }(t)\right\Vert
_{H}^{2}+\delta \int_{0}^{t}\left\Vert (y_{\varepsilon })_{xx}(s)\right\Vert
_{H}^{2}ds\leq \frac{1}{2}\left\Vert \nabla v_{0}\right\Vert _{H}^{2} \\
&&+\int_{0}^{t}(\left\Vert f_{2}(\overline{w}(s))\right\Vert _{H}+\rho
\left\Vert \mbox{sign}_{\varepsilon }(\overline{v}(s)-v^{\ast
}(s))\right\Vert _{H}+\left\Vert f_{1}(\overline{w}(s))y_{\varepsilon
}(s)\right\Vert _{H})\left\Vert (y_{\varepsilon })_{xx}(s)\right\Vert _{H}ds
\\
&\leq &\frac{1}{2}\left\Vert \nabla v_{0}\right\Vert _{H}^{2}+\frac{\delta }{%
2}\int_{0}^{t}\left\Vert (y_{\varepsilon })_{xx}(s)\right\Vert _{H}^{2}ds+%
\frac{3}{2\delta }\left( f_{2M}^{2}LT+\rho
^{2}LT+f_{1M}^{2}\int_0^t\left\Vert y_{\varepsilon }(s)\right\Vert
_{H}^{2}ds\right)
\end{eqnarray*}%
whence%
\begin{eqnarray}
&&\delta \int_{0}^{t}\left\Vert (y_{\varepsilon })_{xx}(s)\right\Vert
_{H}^{2}ds+\left\Vert \nabla y_{\varepsilon }(t)\right\Vert _{H}^{2}
\label{20} \\
&\leq &\left\Vert \nabla v_{0}\right\Vert _{H}^{2}+\frac{3}{\delta }%
(f_{2M}^{2}+\rho ^{2})LT+\frac{3}{2\delta}f_{1M}^{2}\left( \left\Vert
v_{0}\right\Vert _{H}^{2}+2(f_{2M}^{2}+\rho ^{2})LT\right) (e^{\frac{2}{%
\delta_1}{t}}-1)  \nonumber \\
&\leq &\left\{ \left\Vert v_{0}\right\Vert _{V}^{2}\left( 1+\frac{3}{2\delta}%
f_{1M}^{2}\right) +\frac{3}{\delta }(f_{2M}^{2}+\rho ^{2})\left( 1 +
f_{1M}^{2}\right)T L\right\} e^{\frac{2}{\delta_1}T}.  \nonumber
\end{eqnarray}%
The latter together with (\ref{19}) provides%
\begin{equation}
\left\Vert y_{\varepsilon }\right\Vert _{L^{2}(0,T;W)\cap L^{\infty
}(0,T;V)}^{2}\leq \overline{C}^{2}  \label{CR}
\end{equation}
where $\overline{C}^{2}$ is given by 
\begin{eqnarray}
\overline{C}^{2} &=&\left\{ \left( \frac{1}{\delta _{1}}+\frac{1}{\delta _{2}%
}+\frac{3}{2\delta\delta _{2}}f_{1M}^{2}\right) \left\Vert v_{0}\right\Vert
_{V}^{2}\right. ,  \nonumber \\
&&\left. +(f_{2M}^{2}+\rho ^{2})\left( \frac{2}{\delta _{1}}+\frac{3}{\delta
\delta _{2}}\left( 1+f_{1M}^{2}\right) \right) L\right\} e^{\frac{2}{\delta_1%
}T}  \nonumber
\end{eqnarray}%
and $\delta _{2}=\min \{1,\delta \}$. Recalling (\ref{11-0}) we deduce that 
\begin{equation}
\left\Vert y_{\varepsilon }\right\Vert _{L^{\infty }(0,T;V)}\leq R.
\label{20-1}
\end{equation}%
Next, by (\ref{17}) we calculate 
\begin{eqnarray}
&&\left\Vert \frac{dy_{\varepsilon }}{dt}\right\Vert _{L^{2}(0,T;H)}
\label{20-3} \\
&\leq &\left\Vert f_{2}(\overline{w})+\rho \,\mbox{sign}_{\varepsilon }(%
\overline{v}-v^{\ast })-f_{1}(\overline{w})y_{\varepsilon }\right\Vert
_{L^{2}(0,T;H)}+\delta \left\Vert (y_{\varepsilon })_{xx}\right\Vert
_{L^{2}(0,T;H)}  \nonumber \\
&\leq &C(L,T,w_{M},\left\Vert v_{0}\right\Vert _{V},\delta, a,\rho,f_1, f_2)
\nonumber
\end{eqnarray}%
We also recall that, by (\ref{15-0}), $w_{\varepsilon }\in
C^{1}([0,T];C[0,L])$ and 
\begin{equation}
\left\Vert z_{\varepsilon }\right\Vert _{C^{1}[0,T];C[0,L])}\leq
h_{1R}w_{M}+h_{2R}+w_{M}.  \label{20-2}
\end{equation}%
Estimates (\ref{20-1})-(\ref{20-2}) and (\ref{15}) ensure that the solution $%
(y_{\varepsilon },z_{\varepsilon })$ to (\ref{11})-(\ref{14}) belongs to $%
\mathcal{M}.$

Now, let us consider two pairs $(v^{1},w^{1}),$ $(v^{2},w^{2}) \in \mathcal{M%
}$, with the same initial data. We denote by $(y_{\varepsilon
1},z_{\varepsilon 1})$ and $(y_{\varepsilon 2},z_{\varepsilon 2})$ the
corresponding solutions to (\ref{11})-(\ref{14}) and we calculate the
difference of equations (\ref{11}) and (\ref{12}). Namely, we write 
\begin{eqnarray}
&&(y_{\varepsilon 1}-y_{\varepsilon 2})_{t}-\delta (y_{\varepsilon
1}-y_{\varepsilon 2})_{xx}+f_{1}(w^{2})(y_{\varepsilon 1}-y_{\varepsilon 2})
\label{21} \\
&=&f_{2}(w^{1})-f_{2}(w^{2})-\rho (\mbox{sign}_{\varepsilon }(v^{1}-v^{\ast
})-\mbox{sign}_{\varepsilon }(v^{2}-v^{\ast })) -
(f_{1}(w^{1})-f_{1}(w^{2}))y_{\varepsilon1},  \nonumber
\end{eqnarray}%
\begin{equation}
(z_{\varepsilon 1}-z_{\varepsilon 2})_{t}=-(h_{1}(v^{1})z_{\varepsilon
1}-h_{1}(v^{2})z_{\varepsilon 2})+h_{2}(v^{1})-h_{2}(v^{2}),  \label{22}
\end{equation}%
with homogeneous Neumann boundary conditions for $(y_{\varepsilon
1}-y_{\varepsilon 2})$ and zero initial data. Relying on the local Lipschitz
continuity of $f_{i}$ and $h_{i}$, we perform a few calculations in the
right-hand sides of the above equations, denoted $RHS$ and $RHS_{1},$
respectively,%
\[
\left\vert RHS\right\vert \leq \left ( L_{f_{2}}(w_{M}) +
L_{f_{1}}(w_{M})R\right ) \left\vert w^{1}-w^{2}\right\vert +\rho
C_{\varepsilon }\left\vert v^{1}-v^{2}\right\vert 
\]%
and 
\begin{eqnarray*}
\left\vert RHS_{1}\right\vert &\leq &\left\vert
h_{1}(v^{1})-h_{1}(v^{2})\right\vert \left\vert z_{\varepsilon 1}\right\vert
+\left\vert z_{\varepsilon 1}-z_{\varepsilon 2}\right\vert \left\vert
h_{1}(v^{2})\right\vert +\left\vert h_{2}(v^{1})-h_{2}(v^{2})\right\vert \\
&\leq &(L_{h_{1}}(R)w_{M}+L_{h_{2}}(R))\left\vert v^{1}-v^{2}\right\vert
+h_{1R}\left\vert z_{\varepsilon 1}-z_{\varepsilon 2}\right\vert ,
\end{eqnarray*}%
where the constant $C_{\varepsilon }$ depends on $\varepsilon $. We multiply
(\ref{21}) scalarly in $H$ by $(y_{\varepsilon 1}-y_{\varepsilon 2})$ and (%
\ref{22}) by $(z_{\varepsilon 1}-z_{\varepsilon 2})$. We sum up the
resulting equations and then we integrate over $(0,t).$ We get%
\begin{eqnarray*}
&&\frac{1}{2}\left\Vert (y_{\varepsilon 1}-y_{\varepsilon 2})(t)\right\Vert
_{H}^{2}+\delta \int_{0}^{t}\left\Vert \nabla (y_{\varepsilon
1}-y_{\varepsilon 2})(s)\right\Vert _{H}^{2}ds+\frac{1}{2}\left\Vert
(z_{\varepsilon 1}-z_{\varepsilon 2})(t)\right\Vert _{H}^{2} \\
&\leq &\int_{0}^{t}\left\Vert (y_{\varepsilon 1}-y_{\varepsilon
2})(s)\right\Vert _{H}^{2}ds +\int_{0}^{t}\left\Vert (z_{\varepsilon
1}-z_{\varepsilon 2})(s)\right\Vert _{H}^{2}ds \\
&&+\left (\frac{1}{2}L_{f_{2}}^{2}(w_{M}) + L_{f_{1}}^{2}(w_{M})R^2 \right
)\int_{0}^{t}\left\Vert (w^{1}-w^{2})(s)\right\Vert _{H}^{2}ds \\
&&+\left (\frac{1}{2} L_{h_{1}}^{2}(R)w^2_{M} + \frac{1}{2} L_{h_{2}}^{2}(R)
+ \rho ^{2}C_{\varepsilon }^{2}\right )\int_{0}^{t}\left\Vert
(v^{1}-v^{2})(s)\right\Vert _{H}^{2}ds,
\end{eqnarray*}%
so that, defining 
\[
\left\Vert q(t)\right\Vert _{H}^{2}:=\left\Vert (y_{\varepsilon
1}-y_{\varepsilon 2})(t)\right\Vert _{H}^{2}+\left\Vert (z_{\varepsilon
1}-z_{\varepsilon 2})(t)\right\Vert _{H}^{2} 
\]%
by Gronwall's lemma we obtain, for $t\in \lbrack 0,T],$%
\begin{equation}
\left\Vert q(t)\right\Vert _{H}^{2}\leq {\widetilde C_{\varepsilon}}%
\int_{0}^{t}\left( \left\Vert (w^{1}-w^{2})(s)\right\Vert
_{H}^{2}+\left\Vert (v^{1}-v^{2})(s)\right\Vert _{H}^{2}\right) ds,
\label{23}
\end{equation}%
where ${\widetilde C_{\varepsilon}} = C(R,w_{M},\rho,h_1,h_2, f_1, f_2,
\varepsilon)$. In order to show that $\Psi $ is a contraction, we introduce
further the norm $\left\Vert q\right\Vert _{B}=\sup\limits_{t\in \lbrack
0,T]}\left( e^{-\gamma t}\left\Vert q(t)\right\Vert _{H}\right) $ which is
equivalent to the standard norm in $C([0,T];H)$. So that, we multiply (\ref%
{23}) by $e^{-2\gamma t}$ getting 
\begin{eqnarray*}
&&e^{-2\gamma t}\left\Vert q(t)\right\Vert _{H}^{2} \\
&\leq &{\widetilde C_{\varepsilon}}e^{-2\gamma t}\int_{0}^{t}e^{2\gamma
s}e^{-2\gamma s}\left( \left\Vert (w^{1}-w^{2})(s)\right\Vert
_{H}^{2}+\left\Vert (v^{1}-v^{2})(s)\right\Vert _{H}^{2}\right) ds \\
&\leq &{\widetilde C_{\varepsilon}}e^{-2\gamma t}\int_{0}^{t}e^{2\gamma
s}\left( \left\Vert (w^{1}-w^{2})\right\Vert _{B}^{2}+\left\Vert
(v^{1}-v^{2})\right\Vert _{B}^{2}\right) ds \\
&\leq &\frac{{\widetilde C_{\varepsilon}}}{2\gamma }\left( \left\Vert
(w^{1}-w^{2})\right\Vert _{B}^{2}+\left\Vert (v^{1}-v^{2})\right\Vert
_{B}^{2}\right) (1-e^{-2\gamma t}) \\
&\leq &\frac{{\widetilde C_{\varepsilon}}}{2\gamma }\left( \left\Vert
(w^{1}-w^{2})\right\Vert _{B}^{2}+\left\Vert (v^{1}-v^{2})\right\Vert
_{B}^{2}\right) .
\end{eqnarray*}%
Taking the supremum for $t\in \lbrack 0,T],$ and choosing $\gamma $ large
enough such that $2\gamma > {\widetilde C_{\varepsilon}},$ we obtain%
\[
\left\Vert \Psi ({v}^1,{w}^1) - \Psi ({v}^2,{w}^2) \right\Vert
_{B}^{2}=\left\Vert y_{\varepsilon 1}-y_{\varepsilon 2}\right\Vert
_{B}^{2}+\left\Vert z_{\varepsilon 1}-z_{\varepsilon 2}\right\Vert _{B}^{2} 
\]
\[
\leq {\widetilde C_{\varepsilon B}}\left( \left\Vert
(w^{1}-w^{2})(s)\right\Vert _{B}^{2}+\left\Vert (v^{1}-v^{2})(s)\right\Vert
_{B}\right) ,\mbox{ } 
\]%
with ${\widetilde C_{\varepsilon B}}={\widetilde C_{\varepsilon }}/(2\gamma
)<1,$ so that $\Psi $ turns out to be a contraction and to have a unique
fixed point, $\Psi (\overline{v},\overline{w})= (\overline{v},\overline{w})
= (y_{\varepsilon },z_{\varepsilon }).$

This implies that the pair $(v_{\varepsilon},w_{\varepsilon
})=(y_{\varepsilon },z_{\varepsilon })$ is the unique solution to system (%
\ref{1-1})-(\ref{4-1}) satisfying the estimates (\ref{15})-(\ref{15-0}) and (%
\ref{CR})-(\ref{20-2}).

Therefore, along a subsequence (denoted still by $_{\varepsilon })$ we have 
\begin{eqnarray*}
v_{\varepsilon } &\rightarrow &v\mbox{ weakly in }W^{1,2}(0,T;H)\cap
L^{2}(0,T;W), \\
v_{\varepsilon } &\rightarrow &v\mbox{ weak-star in }L^{\infty }(0,T;V), \\
w_{\varepsilon } &\rightarrow &w\mbox{ weak-star in }W^{1,\infty }(0,T;H).
\end{eqnarray*}%
By the Lions-Aubin lemma (see e.g., \cite{Lions-1969}, p. 58) we get%
\begin{eqnarray*}
v_{\varepsilon } &\rightarrow &v\mbox{ strongly in }L^{2}(0,T;V), \\
w_{\varepsilon } &\rightarrow &w\mbox{ strongly in }L^{2}(0,T;H),
\end{eqnarray*}%
implying $v_{\varepsilon }\rightarrow v$ and $w_{\varepsilon }\rightarrow w$
a.e. on $Q.$ By the continuity of $f_{i}$ and $h_{i}$ and by the Lebesgue
dominated convergence theorem we also get 
\[
f_{i}(w_{\varepsilon })\rightarrow f_{i}(w),\quad h_{i}(v_{\varepsilon
})\rightarrow h_{i}(v)\mbox{ strongly in }L^{2}(0,T;H). 
\]%
Also, 
\[
f_{i}(w_{\varepsilon })\rightarrow f_{i}(w),\quad \mbox{ }%
h_{i}(v_{\varepsilon })\rightarrow h_{i}(v)\mbox{ weak-star in }L^{\infty
}(Q). 
\]%
By Arzel\`{a}-Ascoli theorem we still obtain%
\begin{eqnarray*}
v_{\varepsilon }(t) &\rightarrow &v(t)\mbox{ strongly in }H,\mbox{ uniformly
on }[0,T], \\
w_{\varepsilon }(t) &\rightarrow &w(t)\mbox{ strongly in }H,\mbox{ uniformly
on }[0,T].
\end{eqnarray*}%
Also it follows that sign$_{\varepsilon }(v_{\varepsilon }-v^{\ast
})\rightarrow \zeta $ weak-star in $L^{\infty }(Q)$ and since sign is
weakly-strongly closed we get 
\[
\zeta \in \mbox{sign\thinspace }(v-v^{\ast })\mbox{ a.e. }(t,x)\in Q, 
\]%
(see e.g., \cite{vb-springer-2010}, p. 38, Proposition 2.2).

Now, we consider the weak formulation of (\ref{1-1})%
\begin{eqnarray*}
&&\int_{0}^{T}\int_{0}^{L}((v_{\varepsilon })_{t}\psi +\delta \nabla
v_{\varepsilon }\cdot \nabla \psi +f_{1}(w_{\varepsilon })v_{\varepsilon
}\psi +\rho \zeta _{\varepsilon }\psi )dxdt \\
&=&\int_{0}^{T}\int_{0}^{L}f_{2}(w_{\varepsilon })\psi dxdt,\mbox{ for all }%
\psi \in L^{2}(0,T;V),
\end{eqnarray*}%
with $\zeta _{\varepsilon }=$ sign$_{\varepsilon }(v_{\varepsilon }-v^{\ast
})$ and pass to the limit as $\varepsilon $ goes to zero, obtaining (\ref{7}%
). To this end we took into account that 
\begin{eqnarray*}
&&\int_{0}^{T}\int_{0}^{L}(f_{1}(w_{\varepsilon })v_{\varepsilon
}-f_{1}(w)v))\psi dxdt \\
&=&\int_{0}^{T}\int_{0}^{L}(f_{1}(w_{\varepsilon })-f_{1}(w))v_{\varepsilon
}\psi dxdt+\int_{0}^{T}\int_{0}^{L}(v_{\varepsilon }-v)f_{1}(w)\psi
dxdt\rightarrow 0
\end{eqnarray*}%
because $v_{\varepsilon }\rightarrow v$ and $f_{i}(w_{\varepsilon
})\rightarrow f_{i}(w)$ strongly in $L^{2}(0,T;H)$ and $f_{1}(w)\psi \in
L^{2}(0,T;H).$

Passing to the limit in the weak formulation on (\ref{1-1bis}) 
\[
\int_{0}^{T}\int_{0}^{L}((w_{\varepsilon })_{t}\phi +h_{1}(v_{\varepsilon
})w_{\varepsilon }\phi )dxdt=\int_{0}^{T}\int_{0}^{L}h_{2}(v_{\varepsilon
})\phi dxdt,\mbox{ for all }\phi \in L^{2}(0,T;H), 
\]%
and taking into account that $\int_{0}^{T}\int_{0}^{L}(h_{1}(v_{\varepsilon
})w_{\varepsilon }-h_{1}(v)w))\psi dxdt\rightarrow 0$, in a similar way as
above, we get (\ref{7-0}). These two last equations prove that $(v,w)$ is a
solution to (\ref{1})-(\ref{4}).

Moreover, by straightforward calculations using (\ref{1-1bis}) and (\ref{2})
we obtain 
\begin{eqnarray*}
&&\left\vert w_{\varepsilon }-w\right\vert \leq \left\vert \left
(e^{-\int_{0}^{t}h_{1}(v_{\varepsilon}(\sigma ,x))d\sigma } -
e^{-\int_{0}^{t}h_{1}(v)(\sigma ,x)d\sigma }\right) w_{0}(x)\right\vert \\
&& + \int_{0}^{t}\left\vert e^{-\int_{s}^{t}h_{1}(v_{\varepsilon}(\sigma
,x))d\sigma}h_{2}(v_{\varepsilon}(s,x)) - e^{-\int_{s}^{t}h_{1}(v(\sigma
,x))d\sigma} h_{2}(v(s,x))\right\vert ds \\
&&\leq \left ( |w_M|L_{h_1}(R) + T h_2(R)+ L_{h_2}(R)\right )
\int_0^t\left\vert (v_{\varepsilon} - v)(\sigma, x)\right\vert ds.
\end{eqnarray*}%
Then, by integrating on $(0,L)$ we get 
\[
\left\Vert (w_{\varepsilon }-w)(t)\right\Vert _{L^{1}(0,L)}\leq
C(T,R,w_M,h_1,h_2)\left\Vert v_{\varepsilon }-v\right\Vert^2
_{L^{1}(Q)}. 
\]%
Since $w_{\varepsilon }(t)\rightarrow w(t)$ strongly in $H$ uniformly in $t$
we have that $w$ verifies 
\begin{equation}
w(t,x)=e^{-\int_{0}^{t}h_{1}(v(\sigma ,x))d\sigma
}w_{0}(x)+\int_{0}^{t}e^{-\int_{s}^{t}h_{1}(v(\sigma ,x))d\sigma
}h_{2}(v(s,x))ds,  \label{14-2}
\end{equation}%
and it is clear that $w\in C([0,T];C[0,L]),$ since each term is continuous.

For the uniqueness, let $(v_{1},w_{1}),$ $(v_{2},w_{2})$ be two solutions to
(\ref{1})-(\ref{4}) corresponding to the same initial data. We subtract the
equations corresponding to $v_{1}$ and $v_{2},$ 
\begin{eqnarray*}
&&(v_{1}-v_{2})_{t}-\delta (v_{1}-v_{2})_{xx}+\rho (\zeta _{1}-\zeta _{2}) \\
&=&-(f_{1}(w_{1})v_{1}-f_{1}(w_{2})v_{2})+f_{2}(w_{1})-f_{2}(w_{2})
\end{eqnarray*}%
(where $\zeta _{1}\in $ sign\thinspace $(v_{1}-v^{\ast })$ and $\zeta
_{2}\in $ sign\thinspace $(v_{2}-v^{\ast })$ a.e. $(t,x)\in Q)$ and the
equations corresponding to $w_{1}$ and $w_{2},$ 
\[
(w_{1}-w_{2})_{t}=-(h_{1}(v_{1})w_{1}-h_{1}(v_{2})w_{2})+h_{2}(v_{1})-h_{2}(v_{2}). 
\]%
Let us multiply the first difference by $v_{1}-v_{2}$ and the second by $%
w_{1}-w_{2} $, integrate over $(0,t)\times (0,L)$ and sum the resulting
equations. After similar calculations as before, we get 
\begin{eqnarray*}
&&\frac{1}{2}\left\Vert (v_{1}-v_{2})(t)\right\Vert _{H}^{2}+\frac{1}{2}%
\left\Vert (w_{1}-w_{2})(t)\right\Vert _{H}^{2} \\
&\leq &C_{1}\left( \int_{0}^{t}\left\Vert (v_{1}-v_{2})(s)\right\Vert
_{H}^{2}ds+\int_{0}^{t}\left\Vert (w_{1}-w_{2})(s)\right\Vert
_{H}^{2}ds\right)
\end{eqnarray*}%
which yields, by Gronwall's lemma, that $v_{1}(t)=v_{2}(t)$ and $%
w_{1}(t)=w_{2}(t),$ for all $t\in \lbrack 0,T].$ This proves the solution
uniqueness and ends the proof.\hfill $\square $

\medskip

We prove now the occurence of the sliding mode at a finite time $T^{\ast }$.

\medskip

\medskip \noindent \textbf{Theorem 2.2.} \textit{Let} 
\begin{equation}
A=\left\Vert v_{t}^{\ast }\right\Vert _{\infty }+\delta \left\Vert
v_{xx}^{\ast }\right\Vert _{\infty }+f_{1M}\left\Vert v^{\ast }\right\Vert
_{\infty }+f_{2M}+f_{1M}\left\Vert v_{0}-v^{\ast }\right\Vert _{\infty }
\label{25}
\end{equation}%
\textit{and let} 
\begin{equation}
\rho >A+\frac{\left\Vert v_{0}-v^{\ast }\right\Vert _{\infty }}{T}.
\label{26}
\end{equation}%
\textit{Then, for }$T^{\ast }\in \lbrack 0,T]$ \textit{defined as }%
\begin{equation}
T^{\ast }=\frac{\left\Vert v_{0}-v^{\ast }\right\Vert _{\infty }}{\rho -A}
\label{27}
\end{equation}%
\textit{it holds} 
\begin{equation}
v(t,x)=v^{\ast }(t,x), \mbox{ \textit{for all} }t\in \lbrack T^{\ast },T]%
\mbox{ \textit{and all} }x\in \lbrack 0,L].  \label{28}
\end{equation}

\medskip

\noindent \textbf{Proof.} We shall compare the solution to (\ref{1}) with
the solution to the system%
\begin{eqnarray}
q_{t}+\rho \,\mbox{sign\thinspace }q &\ni &A,\mbox{ }t\in (0,T),  \label{29}
\\
q(0) &=&q_{0}=\left\Vert v_{0}-v^{\ast }\right\Vert _{\infty }.
\label{29bis}
\end{eqnarray}%
Since $\frac{A}{\rho }<1$ by (\ref{26}), $\frac{A}{\rho }\in $
sign\thinspace $0$ and one can verify that the solution to (\ref{29}) is 
\textbf{\ }%
\begin{equation}
q(t)=(\left\Vert v_{0}-v^{\ast }\right\Vert _{\infty }-(\rho -A)t)^{+},
\label{30}
\end{equation}%
where $(\cdot )^{+}$ is the positive part. Moreover, it can be noticed that $%
q(T^{\ast })=0$ where $T^{\ast }$ is given by (\ref{27}). Observe that the
function $q$ is positive and decreasing, $\left\vert q(t)\right\vert \leq
\left\vert q(0)\right\vert $ for $t<T^{\ast },$ it reaches the value zero at 
$T^{\ast }$ and remains zero after $T^{\ast }.$ It is clear that due to the
choice (\ref{26}) we have $T^{\ast }<T.$

We denote $p=v-v^{\ast }$ and consider the system%
\begin{equation}
p_{t}-\delta p_{xx}+f_{1}(w)p+\rho \,\mbox{sign\thinspace }p\ni -v_{t}^{\ast
}+\delta v_{xx}^{\ast }-f_{1}(w)v^{\ast }+f_{2}(w),\mbox{ }  \label{31}
\end{equation}%
\begin{equation}
q_{t}-\delta q_{xx}+f_{1}(w)q+\rho \,\mbox{sign\thinspace }q\ni A+f_{1}(w)q
\label{32}
\end{equation}%
\begin{equation}
p(0)=v_{0}-v^{\ast },\mbox{ }q(0)=\left\Vert v_{0}-v^{\ast }\right\Vert
_{\infty }  \label{33}
\end{equation}%
with homogeneous Neumann boundary conditions both for $p$ and $q.$ Observe
that since $q$ depends only on time (cf. (\ref{30})), then (\ref{29})-(\ref%
{29bis}) is equivalent to (\ref{32})-(\ref{33}-ii).

We subtract (\ref{32}) from (\ref{31}) and multiply the difference equation
scalarly in $H$ by the positive part $(p-q)^{+}$ and integrate over $(0,t).$
By few calculations and majorating the right-hand side of the difference
equation, we obtain 
\begin{eqnarray*}
&&\frac{1}{2}\left\Vert (p-q)^{+}(t)\right\Vert _{H}^{2}+\delta
\int_{0}^{t}\left\Vert \nabla (p-q)^{+}(s)\right\Vert
_{H}^{2}ds+\int_{0}^{t}\int_{0}^{L}f_{1}(w)\left( (p-q)^{+}(s)\right) ^{2}ds
\\
&&+\rho \int_{0}^{t}\int_{0}^{L}(\zeta _{p}-\zeta _{q})(p-q)^{+}dxds \\
&= &\int_{0}^{t}\int_{0}^{L}\left( -v_{t}^{\ast }+\delta v_{xx}^{\ast
}-f_{1}(w)v^{\ast }+f_{2}(w)-f_{1}(w)q-A\right) (p-q)^{+}dxds \\
&\leq &\int_{0}^{t}\int_{0}^{L}(\left\Vert v_{t}^{\ast }\right\Vert _{\infty
}+\delta \left\Vert v_{xx}^{\ast }\right\Vert _{\infty }+f_{1M}\left\Vert
v^{\ast }\right\Vert _{\infty }+f_{2M}+f_{1M}\left\Vert q_{0}\right\Vert
_{\infty }-A)(p-q)^{+}dxds=0,
\end{eqnarray*}%
by (\ref{25}), where $\zeta _{p}\in $ sign$\,p$ and $\zeta _{q}\in $ sign$%
\,q.$ We took into account that $p(0)-q(0)=v_{0}-v^{\ast }-\left\Vert
v_{0}-v^{\ast }\right\Vert _{\infty }\leq 0$ and so $(p(0)-q(0))^{+}=0.$
From here it follows that $p(t)\leq q(t)$ for all $t\in \lbrack 0,T].$

Now, we add the equations for $p$ and $q$ and multiply their sum by $%
-(p+z)^{-}$ and integrate over $(0,t)\times (0,L).$ Taking into account that
sign$z=-$sign$(-z)$ we obtain%
\begin{eqnarray*}
&&\frac{1}{2}\left\Vert (p+q)^{-}(t)\right\Vert _{H}^{2}+\delta
\int_{0}^{t}\left\Vert \nabla (p+q)^{-}(s)\right\Vert
_{H}^{2}ds+\int_{0}^{t}\int_{0}^{L}f_{1}(w)\left( (p+q)^{-}(s)\right) ^{2}ds
\\
&&-\rho \int_{0}^{t}\int_{0}^{L}(\zeta _{p}-\zeta _{q})(p-(-q))^{-}dxds \\
&=&-\int_{0}^{t}\int_{0}^{L}\left( -v_{t}^{\ast }+\delta v_{xx}^{\ast
}-f_{1}(w)v^{\ast }+f_{2}(w)+f_{1}(w)q+A\right) (p-q)^{-}dxds,
\end{eqnarray*}%
where $\zeta _{q}\in $ sign$\,(-q)$. Observe that ($p(0)+q(0))=
v_{0}-v^{\ast }+\left\Vert v_{0}-v^{\ast }\right\Vert _{\infty } \geq 0$ and
so $(p(0)+q(0))^{-}=0.$ Thus, 
\begin{eqnarray*}
&&\left\Vert (p+q)^{-}(t)\right\Vert _{H}^{2}\leq
\int_{0}^{t}\int_{0}^{L}\left( v_{t}^{\ast }-\delta v_{xx}^{\ast
}+f_{1}(w)v^{\ast }-f_{2}(w)-f_{1}(w)q-A\right) (p+q)^{-}dxds \\
&\leq &\int_{0}^{t}\int_{0}^{L}(\left\Vert v_{t}^{\ast }\right\Vert _{\infty
}+\delta \left\Vert v_{xx}^{\ast }\right\Vert _{\infty }+f_{1M}\left\Vert
v^{\ast }\right\Vert _{\infty }+f_{2M}+f_{1M}\left\Vert q_{0}\right\Vert
_{\infty }-A)(p+q)^{-}dxds=0,
\end{eqnarray*}%
implying that $p(t)\geq -q(t)$ for all $t\in \lbrack 0,T].$

Finally we have obtained that $\left\vert p(t)\right\vert =\left\vert
v(t)-v^{\ast }(t)\right\vert \leq q(t),$ so that $v(t)-v^{\ast }(t)=0$ for $%
t\geq T^{\ast },$ which yields (\ref{28}), as claimed. Finally, we observe
that $A$ is not a sharp value (it could be smaller) but here the objective
was to prove its existence.\hfill $\square $

\medskip

\section{The complete system}

\setcounter{equation}{0}

Relying on the results previously obtained we can pass to the complete
Hodgkin-Huxley system (\ref{HH-new-1})-(\ref{ci}) and assume:

\begin{itemize}
\item[(i)$_{1}$] the functions $f_{i}$ and $h_{i}^{k},$ $i=1,2,$ $k=n,$ $m,$ 
$h,$ are locally Lipschitz continuous, that is, for any $M$ positive, and
for any $\,r,$ $r_{1}$, $r_{2},$ $r_{3},$ $\overline{r},$ $\overline{r}_{1},$
$\overline{r}_{2},$ $\overline{r}_{3}\in \mathbb{R}$, $\left\vert
r_{i}\right\vert \leq M,$ $\left\vert \overline{r}_{i}\right\vert \leq M,$
there exist $L_{f_{i}}(M)$ and $L_{h_{i}^{k}}(M)$ positive, such that 
\begin{eqnarray}
\left\vert f_{i}(r_{1},r_{2},r_{3})-f_{i}(\overline{r}_{1},\overline{r}_{2},%
\overline{r}_{3})\right\vert &\leq
&L_{f_{i}}(M)\sum\limits_{j=1}^{3}\left\vert r_{j}-\overline{r}%
_{j}\right\vert ,  \label{34} \\
\left\vert h_{i}^{k}(r)-h_{i}^{k}(\overline{r})\right\vert &\leq
&L_{h_{i}^{k}}(M)\left\vert r-\overline{r}\right\vert ,\mbox{ }i=1,2; 
\nonumber
\end{eqnarray}

\item[(ii)$_{1}$] there exists $a>0$ such that 
\begin{equation}
0<a\leq f_{1}(r),\mbox{ }0<h_{i}^{k}(r),\mbox{ }i=1,2,\mbox{ for all }r\in 
\mathbb{R};\mbox{ }k=n,\mbox{ }m,\mbox{ }h;  \label{35}
\end{equation}

\item[(iii)$_{1}$] 
\begin{equation}
n_{M}:=\sup_{r\in \mathbb{R}}\frac{h_{2}^{n}(r)}{h_{1}^{n}(r)},\mbox{ }%
m_{M}:=\sup_{r\in \mathbb{R}}\frac{h_{2}^{m}(r)}{h_{1}^{m}(r)},\mbox{ }%
h_{M}:=\sup_{r\in \mathbb{R}}\frac{h_{2}^{h}(r)}{h_{1}^{h}(r)};  \label{38}
\end{equation}%
\begin{equation}
v_{0}\in V,\mbox{ }n_{0},\mbox{ }m_{0},\mbox{ }h_{0}\in C[0,L],  \label{36}
\end{equation}%
\begin{equation}
n_{0}\in \lbrack 0,n_{M}],\mbox{ }m_{0}\in \lbrack 0,n_{M}],\mbox{ }h_{0}\in
\lbrack 0,n_{M}];  \label{37}
\end{equation}

\item[(iv)$_{1}$] 
\begin{equation}
v^{\ast }\in W^{1,\infty }(0,T;H)\cap L^{\infty }(0,T;W).
\end{equation}
\end{itemize}

\medskip

\noindent \textbf{Definition 3.1.} We call a solution to system (\ref%
{HH-new-1})-(\ref{ci}) a vector $(v,n,m,h)$%
\begin{eqnarray*}
v &\in &C([0,T];H)\cap L^{2}(0,T;V)\cap W^{1,2}(0,T;H)\cap L^{\infty }(Q), \\
n,\mbox{ }m,\mbox{ }h &\in &L^{\infty }(Q)\cap W^{1,2}(0,T;H)\mbox{, }
\end{eqnarray*}%
which satisfies 
\begin{eqnarray*}
&&\int_{0}^{T}\int_{0}^{L}(v_{t}\psi +\delta \nabla v\cdot \nabla \psi
+f_{1}(n,m,h)v\psi +\rho \zeta \psi )dxdt \\
&=&\int_{0}^{T}\int_{0}^{L}f_{2}(n,m,h)\psi dxdt,\mbox{ for all }\psi \in
L^{2}(0,T;V),\mbox{ } \\
\zeta (t,x) &\in &\mbox{sign}\,(v(t,x)-v^{\ast }(t,x)),\mbox{ a.e. }(t,x)\in
(0,T)\times (0,L),
\end{eqnarray*}%
and%
\[
\int_{0}^{T}\int_{0}^{L}(n_{t}\phi +h_{1}^{n}(v)n\phi
)dxdt=\int_{0}^{T}\int_{0}^{L}h_{2}^{n}(v)\phi dxdt,\mbox{ for all }\phi \in
L^{2}(0,T;H), 
\]%
\[
\int_{0}^{T}\int_{0}^{L}(m_{t}\phi +h_{1}^{m}(v)m\phi
)dxdt=\int_{0}^{T}\int_{0}^{L}h_{2}^{m}(v)\phi dxdt,\mbox{ for all }\phi \in
L^{2}(0,T;H), 
\]%
\[
\int_{0}^{T}\int_{0}^{L}(h_{t}\phi +h_{1}^{h}(v)h\phi
)dxdt=\int_{0}^{T}\int_{0}^{L}h_{2}^{h}(v)\phi dxdt,\mbox{ for all }\phi \in
L^{2}(0,T;H), 
\]%
together with the initial conditions (\ref{ci}).

\medskip

\noindent \textbf{Theorem 3.1.} \textit{Let} (1)$_{1}$-(iv)$_{1}$\textit{\
hold. Then, problem} (\ref{HH-new-1})-(\ref{ci}) \textit{has a unique
solution, which has the supplementary regularity} 
\[
v\in L^{\infty }(0,T;V)\cap L^{2}(0,T;W),\mbox{ }n,\mbox{ }m,\mbox{ }h\in
C([0,T];C[0,L])\cap W^{1,\infty }(0,T;H). 
\]%
\textit{Moreover, if} 
\[
\rho >A+\frac{\left\Vert v_{0}-v^{\ast }\right\Vert _{\infty }}{T}. 
\]%
\textit{with }$A$ \textit{as in} (\ref{25}), \textit{then for } $T^{\ast }
\in (0,T)$ defined as%
\[
T^{\ast }=\frac{\left\Vert v_{0}-v^{\ast }\right\Vert _{\infty }}{\rho -A} 
\]%
\textit{it holds} 
\[
v(t,x)=v^{\ast }(t,x)\mbox{ \textit{for all} }t\in \lbrack T^{\ast },T]\mbox{
\textit{and all} }x\in \lbrack 0,L]. 
\]

\medskip

\noindent \textbf{Proof. }In\textbf{\ }(\ref{HH-new-1})-(\ref{ci}) we set 
\[
w=\left( 
\begin{array}{c}
n \\ 
m \\ 
h%
\end{array}%
\right) ,\mbox{ }h_{1}(v)=\left( 
\begin{array}{ccc}
h_{1}^{n}(v)n & 0 & 0 \\ 
0 & h_{1}^{m}(v) & 0 \\ 
0 & 0 & h_{1}^{h}(v)%
\end{array}%
\right) ,\mbox{ }h_{2}(v)=\left( 
\begin{array}{c}
h_{2}^{n}(v) \\ 
h_{2}^{m}(v) \\ 
h_{2}^{h}(v)%
\end{array}%
\right) 
\]%
and so system (\ref{HH-new-1})-(\ref{ci}) can be written in the form (\ref{1}%
)-(\ref{4}). Also, we observe that if $\overline{w}=\left( 
\begin{array}{c}
\overline{n} \\ 
\overline{m} \\ 
\overline{h}%
\end{array}%
\right) $ and $\left\vert n\right\vert \leq M,$ $\left\vert \overline{n}%
\right\vert \leq M,$ $\left\vert m\right\vert \leq M,$ $\left\vert \overline{%
m}\right\vert \leq M,$ $\left\vert h\right\vert \leq M,$ $\left\vert 
\overline{h}\right\vert \leq M$, then 
\begin{eqnarray*}
\left\vert w\right\vert \leq \sqrt{3}M, \left\vert \overline{w}\right\vert
&\leq& \sqrt{3}M, \\
\left\vert f_{1}(w)-f_{1}(\overline{w})\right\vert &\leq
&4M^{3}g_{K}\left\vert n-\overline{n}\right\vert +3M^{2}Mg_{Na}\left\vert m-%
\overline{m}\right\vert +M^{3}g_{Na}\left\vert h-\overline{h}\right\vert \\
&\leq &L_{f_{1}}\left\vert w-\overline{w}\right\vert ,\mbox{ }L_{f_{1}}=\max
\{4M^{3}g_{K},3M^{2}Mg_{Na},M^{3}g_{Na}\}, \\
\left\vert f_{1}(w)-f_{1}(\overline{w})\right\vert &\leq
&L_{f_{2}}\left\vert w-\overline{w}\right\vert ,\mbox{ }L_{f_{2}}=\max
\{4M^{3}g_{K}V_{K},3M^{2}Mg_{Na}V_{Na},M^{3}g_{Na}V_{Na}\}.
\end{eqnarray*}%
Moreover, $h_{2}(v)$ and each column vector in $h_{1}(v)$ are locally
Lipschitz due to the same properties of $h_{i}^{k}(v)$ by (\ref{34}). Then,
one can apply Theorem 2.1 and take as set $\mathcal{M}$ the following%
\begin{eqnarray*}
\mathcal{M} &=&\left\{ (v,n,m,h)\in \left ( C([0,T];H) \cap L^\infty(0,T;V)
\right ) \times (C([0,T];H) \cap L^\infty(Q))^{3};\right. \\
&&\left. \left\Vert v\right\Vert _{L^{\infty }(0,T;V)}\leq R,\mbox{ }%
\left\Vert n\right\Vert _{L^{\infty }(Q)}\leq n_{M},\mbox{ }\left\Vert
m\right\Vert _{L^{\infty }(Q)}\leq m_{M},\mbox{ }\left\Vert h\right\Vert
_{L^{\infty }(Q)}\leq h_{M}\right. \}.
\end{eqnarray*}%
Here, we set $f_{iM}:=f_{i}(0,0,0)+L_{f_{i}}\left( n_{M}+m_{M}+h_{M}\right)
. $ Then, Theorem 2.2 can be applied to get the result. \hfill $\square $

\section{Numerical simulations}

\setcounter{equation}{0}

We present some numerical simulations intended to show the feature of the HH
system evolution controlled by the relay controller and to put into evidence
the sliding mode behavior.

The numerical simulations have been done for the complete system (\ref%
{HH-new-1})-(\ref{ci}) with $I_{C}(t,x)=\rho \,$sign$_{\varepsilon
}(v(t,x)-v^{\ast }(t,x)),$ which was solved by an interactive technique.
Thus, the numerical solution is computed for the approximating system, but
for simplicity, we shall refer later to $v,$ $n,$ $m,$ $h$ without the
subscript $\varepsilon .$

We considered the domain $Q=[0,T]\times \lbrack 0,L]$ with $L=1$ and $T\in
\{100,200,400\}$ and the approximation of the multivalued function sign
given by 
\[
\mbox{sign}_{\varepsilon }r=\left\{ 
\begin{array}{l}
1,\mbox{ \ \ \ }r>\varepsilon \\ 
\frac{r}{\varepsilon },\mbox{ \ \ \ }r\in (-\varepsilon ,\varepsilon ) \\ 
-1,\mbox{ \ }r<-\varepsilon%
\end{array}%
\right. 
\]%
with $\varepsilon =10^{-4}.$

To solve the ordinary differential equations, ODE, in $n$, $m$ and $h$ with
the corresponding initial conditions we used the \textit{ode45} solver from
Matlab. The solver is based on an explicit Runge-Kutta (4,5) formula, the
Dormand-Prince pair (\cite{Dormand}). It is a one-step solver for computing\
a value $y(t_{n})$ and it needs only the solution at the immediately
preceding time point $y(t_{n-1}).$

The numerical solution of the initial-homogeneous Neumann boundary value
problem for the 1-D parabolic equation in $v$ was computed with the\textit{\
pdepe} Matlab solver. The solver discretizes the space using a given $%
xmesh=[x_{1}<...<x_{maxX}]$ and integrates the resulting ODE to obtain
approximate solutions at times specified by a vector of points $%
tmesh=[t_{1}<...<t_{maxT}]$ for all points in $xmesh$ (see \cite{Skeel}).
The time integration is done with \textit{ode15s}, a variable order Matlab
solver based on the numerical differentiation formulas (see \cite{Shampine}).

To discretize the space interval $[0,L]$ we considered $maxX$ steps, with $%
maxX\in \{5,25\}$, and to discretize the time interval $[0,T]$ we took $maxT$
steps, with $maxT\in \{20,40,200\}$.

We used the following stopping criteria of the algorithm:

\[
\Vert v^{iter}-v^{iter-1}\Vert _{\infty }<10^{-3}, 
\]%
where $v^{iter}$ is the value of $v$ at iteration $iter$;

\[
iter>Niter, 
\]%
that is, the number of iterations exceeds $Niter(=100)$;

\[
eTime>NTime, 
\]%
that is, the elapsed time exceeds the limit $NTime(=900CPU)$.

The idea of the algorithm is to solve iteratively, on blocks, the system $%
(v,n,m,h)$ untill one of the stopping conditions is fulfilled. In our case
we consider two blocks, the first containing the PDE in $v$, the second
block containing the ODE system $(n,m,h)$.

A solving iteration consists in: computing a numerical solution of $v$ using
the previous values of $n$, $m$, $h$ and after, to plug-in the obtained
value of $v$ in the right terms of ODE for $n$, $m$, $h$ and solve the
equations from the second block.

The pseudo-code of algorithm is:

\textbf{Step 1.} Discretize the domain $Q$, compute the initial and boundary
conditions and initialize the solving loop.

- Generate the vectors $tmesh$, $xmesh$ and the matrix $meshgrid=%
\{(t_{i},x_{j})\}_{i=1,.,maxT}^{j=1,.,maxX}$ using the given $maxT$ and $%
maxX $.

- Evaluate, in the grid points $meshgrid$ the functions $v_{0},$ $n_{0},$ $%
m_{0},$ $h_{0}$ and $v^{\ast }$.

- $iter=0.$

- $n^{iter}=n_{0},$ $m^{iter}=m_{0},$ $h^{iter}=h_{0},$ $v^{iter}=v_{0},$ $%
eTime=0.$

\textbf{Step 2.}

- Compute $f_{1}(n^{iter},m^{iter},h^{iter})$ and $%
f_{1}(n^{iter},m^{iter},h^{iter}).$

\textbf{Step 3.} Iteration loop.

- $iter=iter+1$.

- Compute $v^{iter}$ using the \textit{pdepe} solver, with the initial datum 
$v_{0}$ and the boundary data.

\textbf{Step 4.}

- Evaluate $h_{n}^{1}(v^{iter}),$ $h_{m}^{1}(v^{iter}),h_{h}^{1}(v^{iter})$
and $h_{n}^{2}(v^{iter}),$ $h_{m}^{2}(v^{iter}),h_{h}^{2}(v^{iter}).$

\textbf{Step 5.} For each point in $xmesh$, compute:

- $n^{iter}$ using the \textit{ode45} solver, with the initial datum $n_0$, for all points in $tmesh$

- $m^{iter}$ using the \textit{ode45} solver, with the initial datum $m_0$, for all points in $tmesh$

- $h^{iter}$ using the \textit{ode45} solver, with the initial datum $h_0$,  for all points in $tmesh.$

\textbf{Step 6.} Check the stopping criteria.

If one of the three conditions is met the algorithm is finished, otherwise

- $eTime$ = elapsed time from the beginning of the iteration loop.

- Go to Step 2.

The algorithm converges, i.e., the first stopping condition is met, in most
cases in a short time (dozen of seconds). The second and third stopping
criterion is used in the non-stabilization cases.

The initial condition was selected by assuming that the sodium channel
inactivation ratio is higher than that of the activation, 
\[
n_{0}=0.45,\mbox{ }m_{0}=0.03,\mbox{ }h_{0}=0.397. 
\]%
For the initial $v_{0}$ various values were considered and they are
indicated in the figures.

Most part of the parameters used in the computations is the same as in \cite%
{Doruk-2010}:%
\begin{eqnarray}
g_{K} &=&36,\mbox{ }g_{Na}=120,\mbox{ }g_{l}=0.3,  \label{val} \\
V_{K} &=&-12,\mbox{ }V_{Na}=115,\mbox{ }V_{l}=10.613  \nonumber \\
\delta  &=&0.1,\mbox{ }C_{M}=0.91,\mbox{ }\rho \geq 0.  \nonumber
\end{eqnarray}%
In some figures the values of $\rho ,$ $g_{K}$ and $T$ differ from those
before and they are specified in the captions.

The graphics plotted in all figures show (from left to write) the time
evolution of $v$ at specified fixed $x,$ the surface $v(t,x)$ and the time
evolution of the proportions of the activating molecules of the potassium $n$
$(K)$ and natrium $m$ $(Na)$ channels and the proportion of the inactivating
molecules of molecules of sodium $h$ $(Iso)$ at specified fixed $x.$

    \begin{center}
    \begin{tabular}{|c|c|c|}
    \hline
    \includegraphics[width=2.0in]{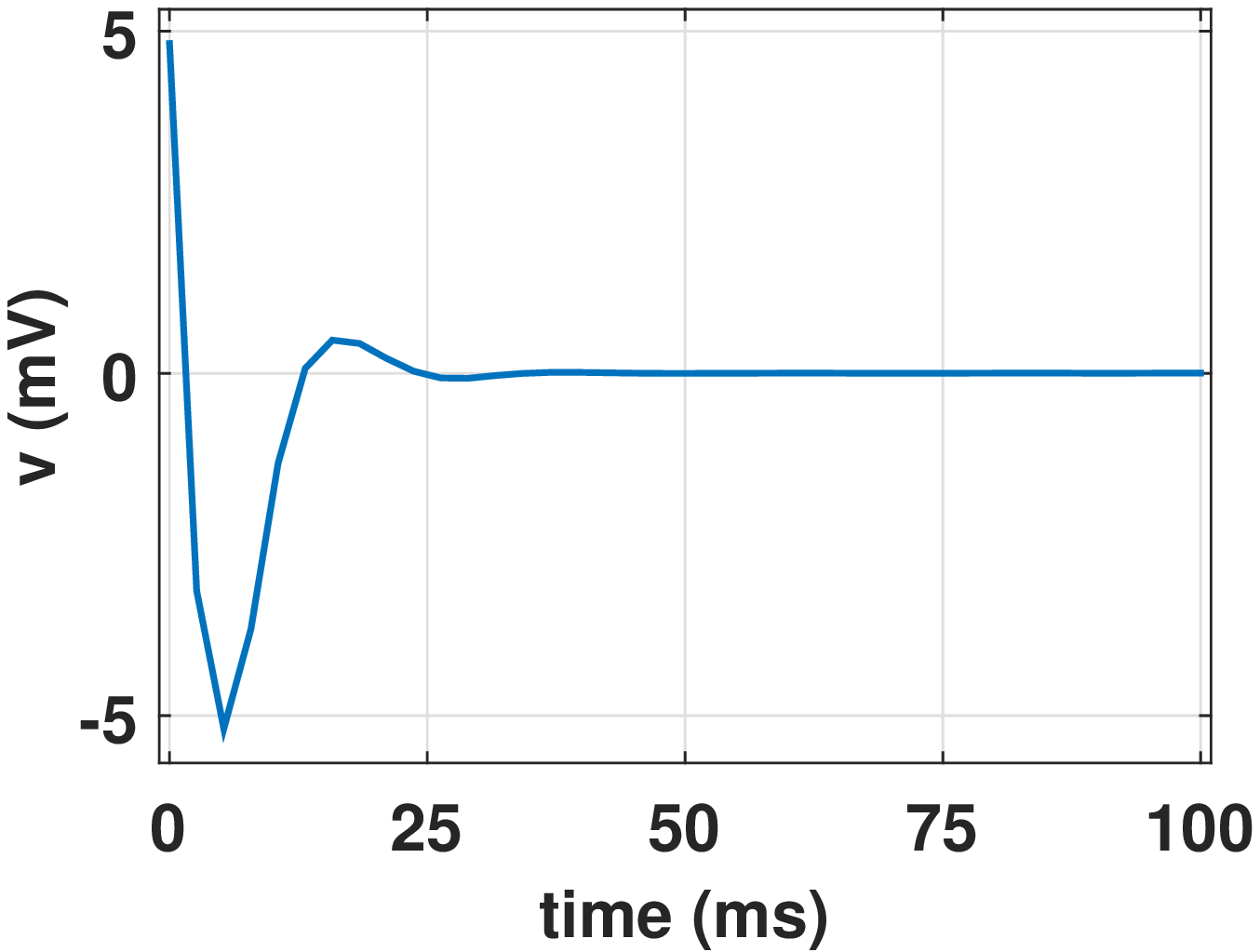} & \includegraphics[width=2.0in]{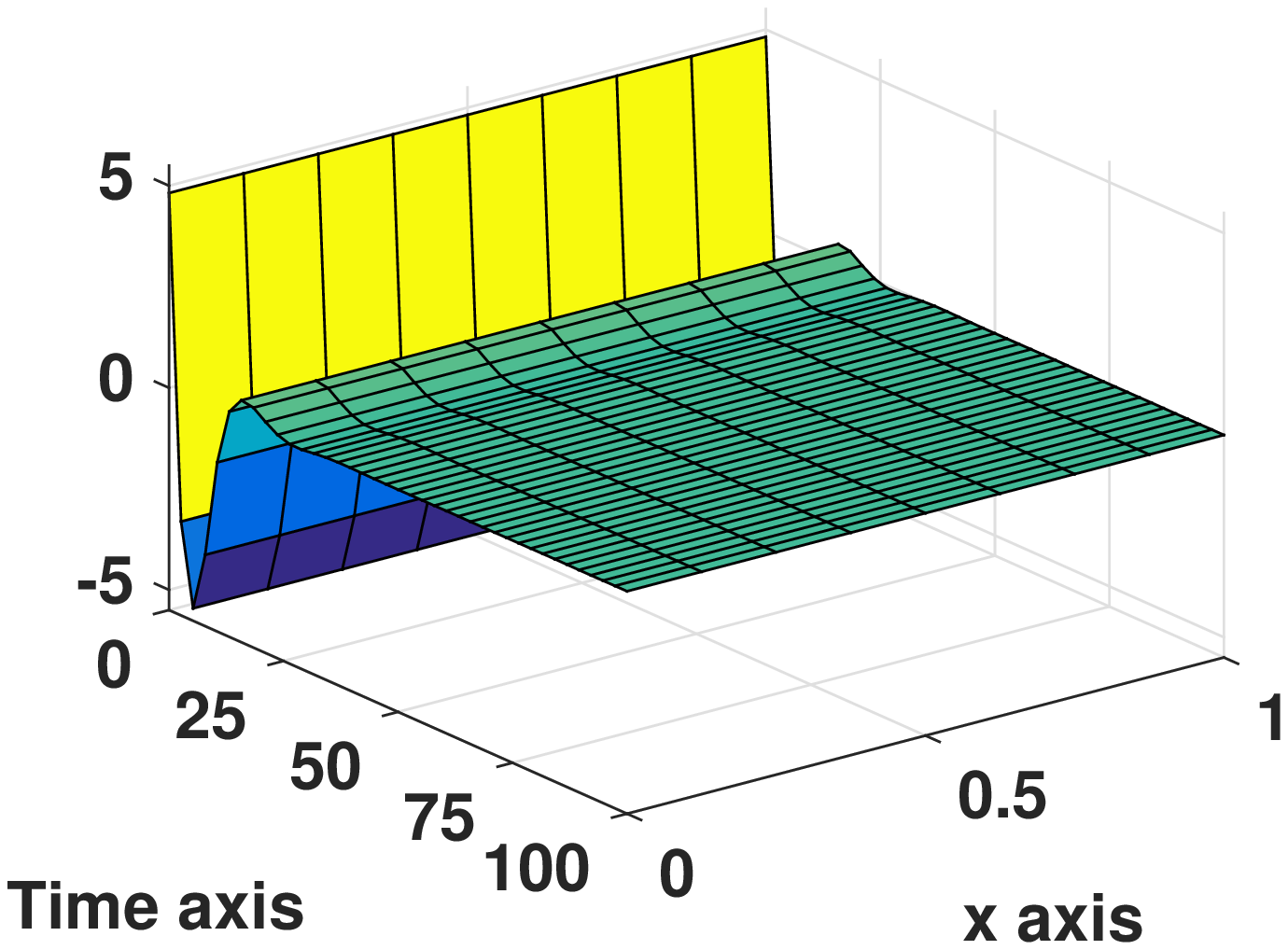} & \includegraphics[width=2.0in]{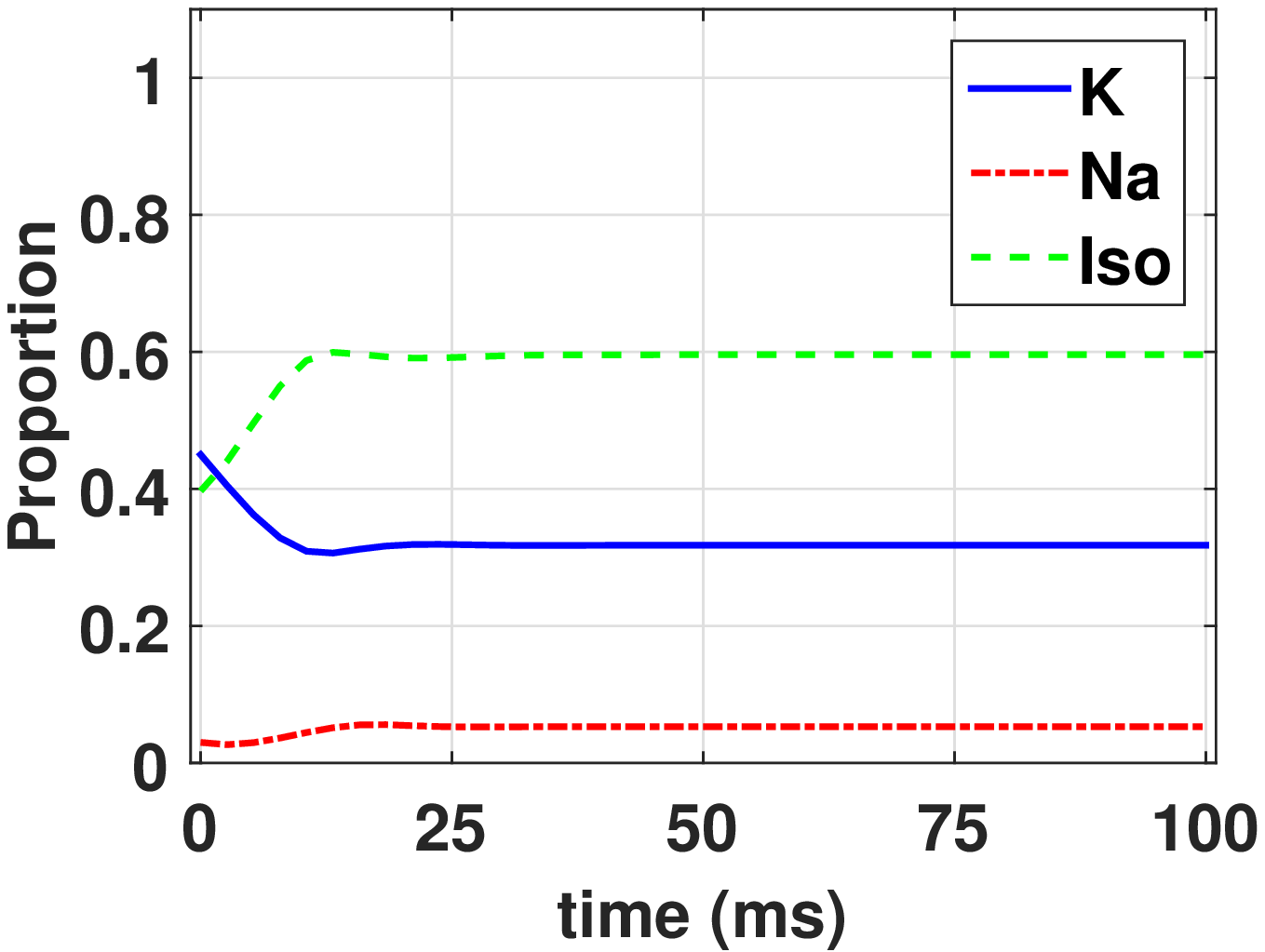} \\
        \hline
    \end{tabular}\\

    {\bf{Fig. 1.}} Graphics $v(t,0)$ (left), $v(t,x)$ (center), $n$, $m$, $h$ (right) for $v_0 = 4.82$, $v^* = 0$, $\rho = 0$
    \end{center}

The values (\ref{val}) are the values of membrane channel conductance which,
for $\rho =0,$ do not lead to an unstable membrane potential response (see 
\cite{Doruk-2010}). This situation is illustrated in Fig. 1. However, for
the same values, the computations for $\rho =20$ show in Fig. 2 a quicker
stabilization.

 	\begin{center}
    \begin{tabular}{|c|c|c|}
    \hline
    \includegraphics[width=2.0in]{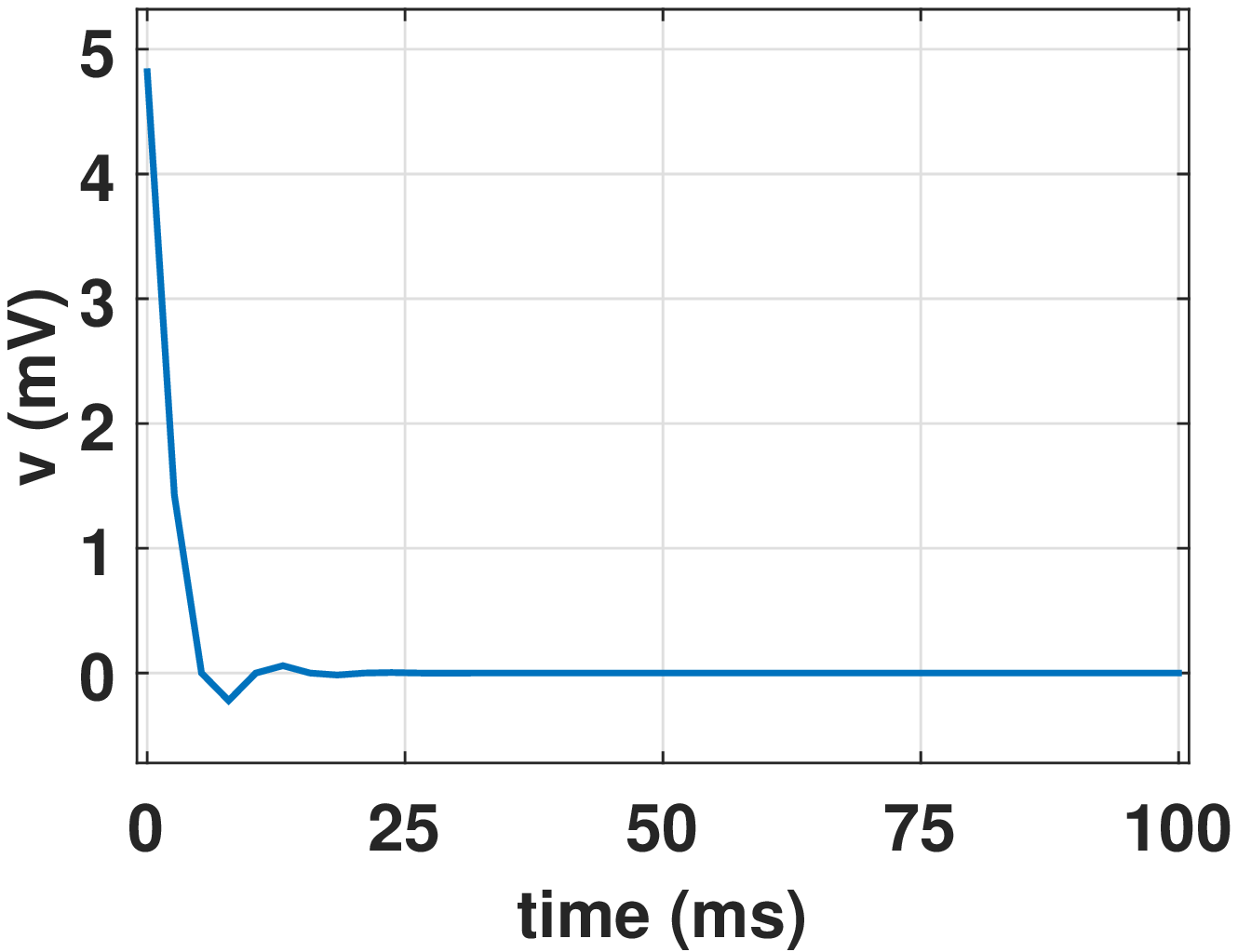} & \includegraphics[width=2.0in]{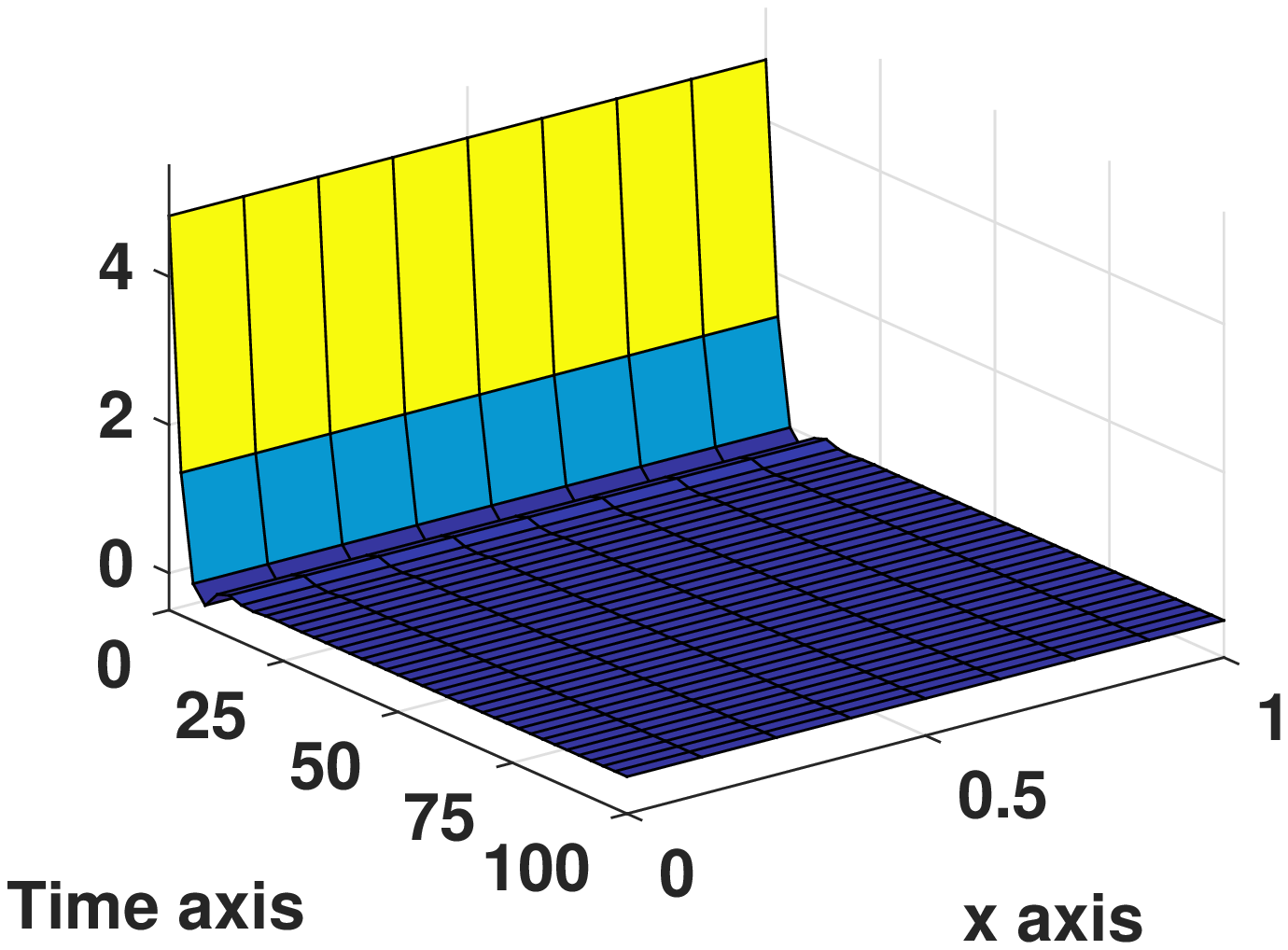} & \includegraphics[width=2.0in]{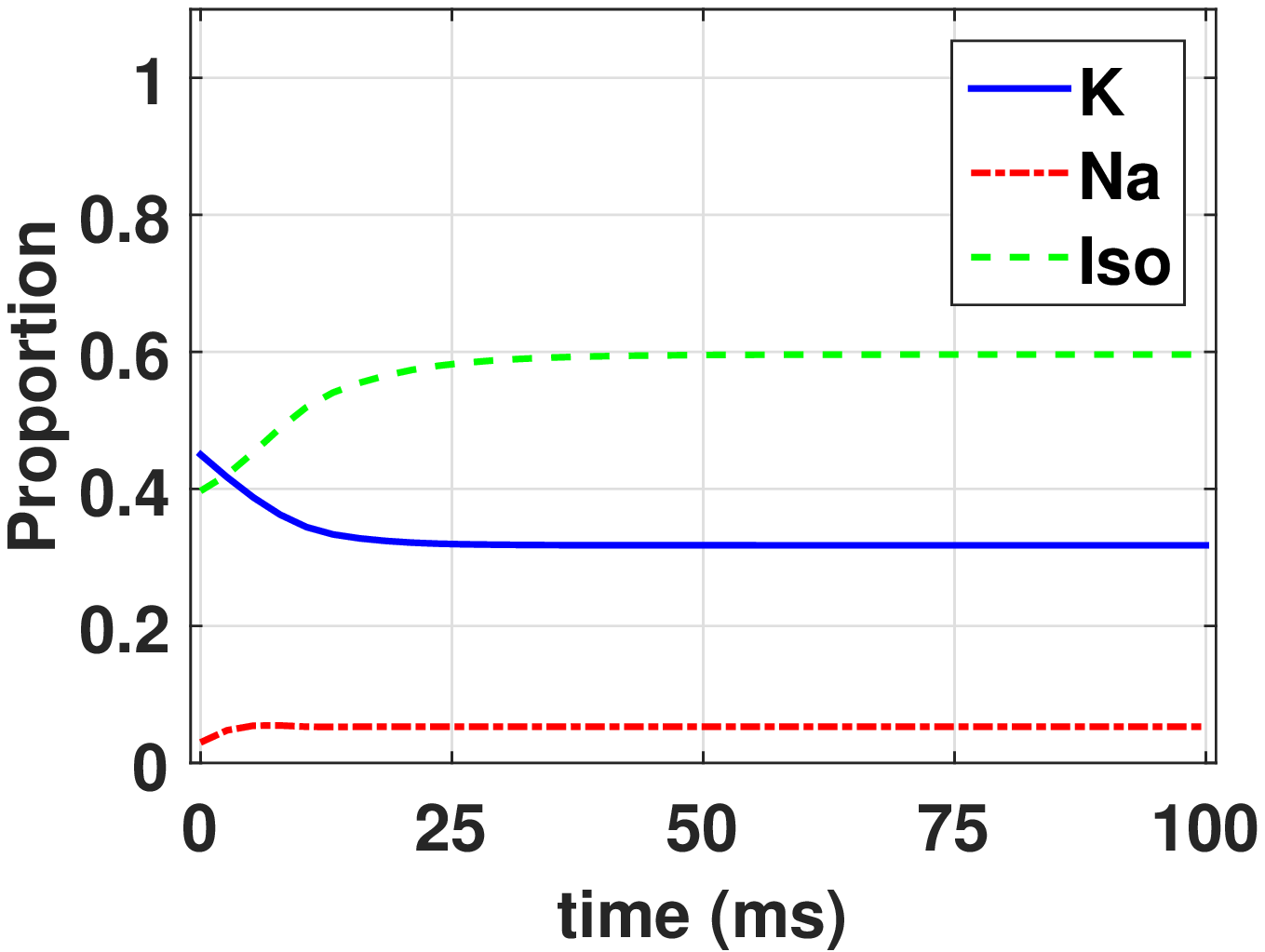} \\
        \hline
    \end{tabular}\\

    {\bf{Fig. 2.}} Graphics $v(t,0)$ (left), $v(t,x)$ (center), $n$, $m$, $h$ (right) for $v_0 = 4.82$, $v^* = 0$, $\rho = 20$
    \end{center}

In order to illustrate the theory that allows the solution to reach a
periodic sliding mode we present Fig. 3 which describes such a situation for 
$v^{\ast }(t,x)=0.5\sin \left( \frac{4}{\pi }t\right) +0.6$ and the same
values (\ref{val}).

For some deviation in conductance parameters, the situation can change with
respect to the first case (see \cite{Doruk-2010}). Thus, for the potassium
channel conductance $g_{K}=3.8229,$ the stabilization does not occur and as
a matter of fact one observes in Fig. 4 a late firing behavior. The desired
stabilization is obtained for a suitable $\rho =20$ and it is illustrated in
Fig. 5. All the other parameters are those from (\ref{val}).
\medskip

	\begin{center}
    \begin{tabular}{|c|c|c|}
    \hline
    \includegraphics[width=2.0in]{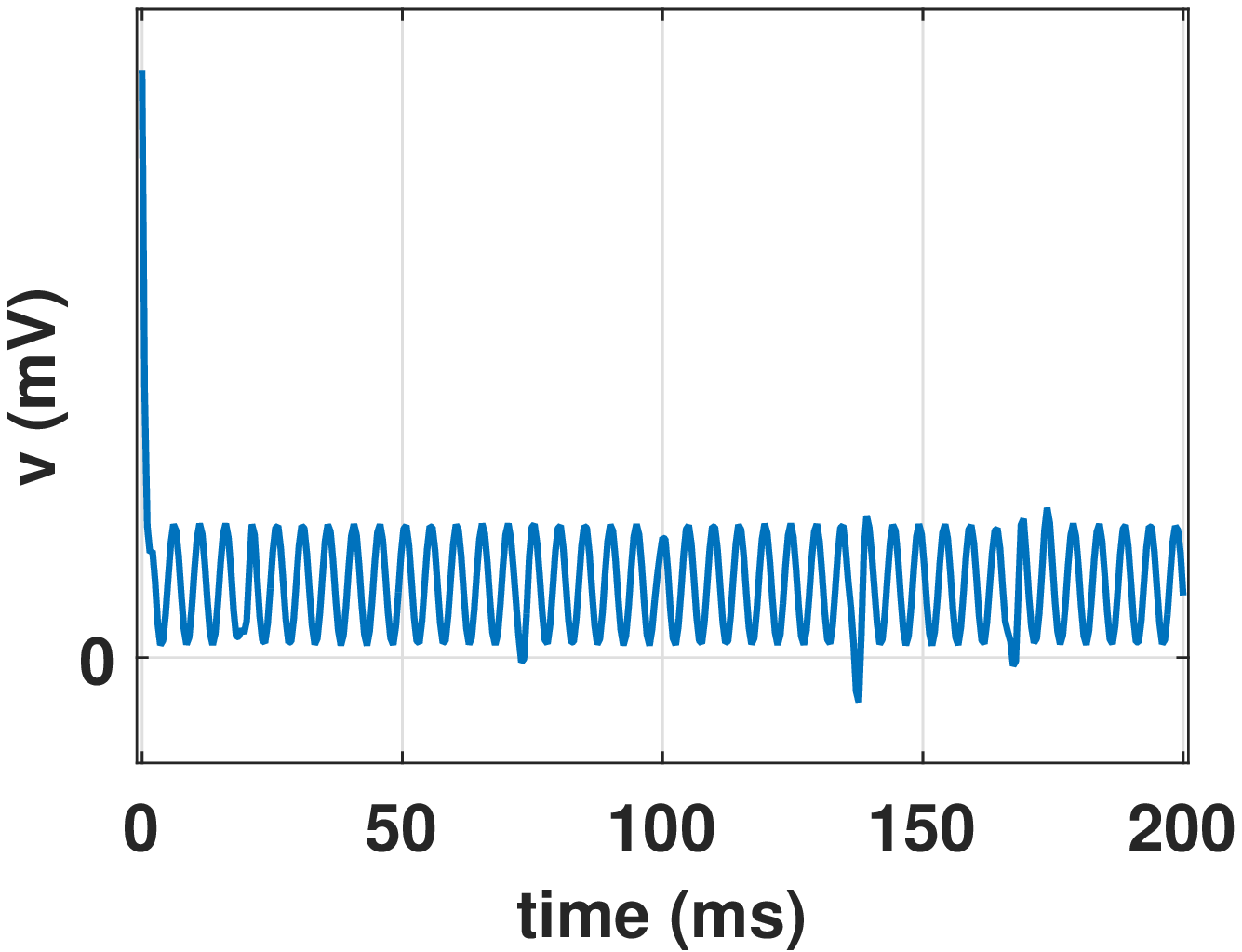} & \includegraphics[width=2.0in]{IIa_v-3D.eps} & \includegraphics[width=2.0in]{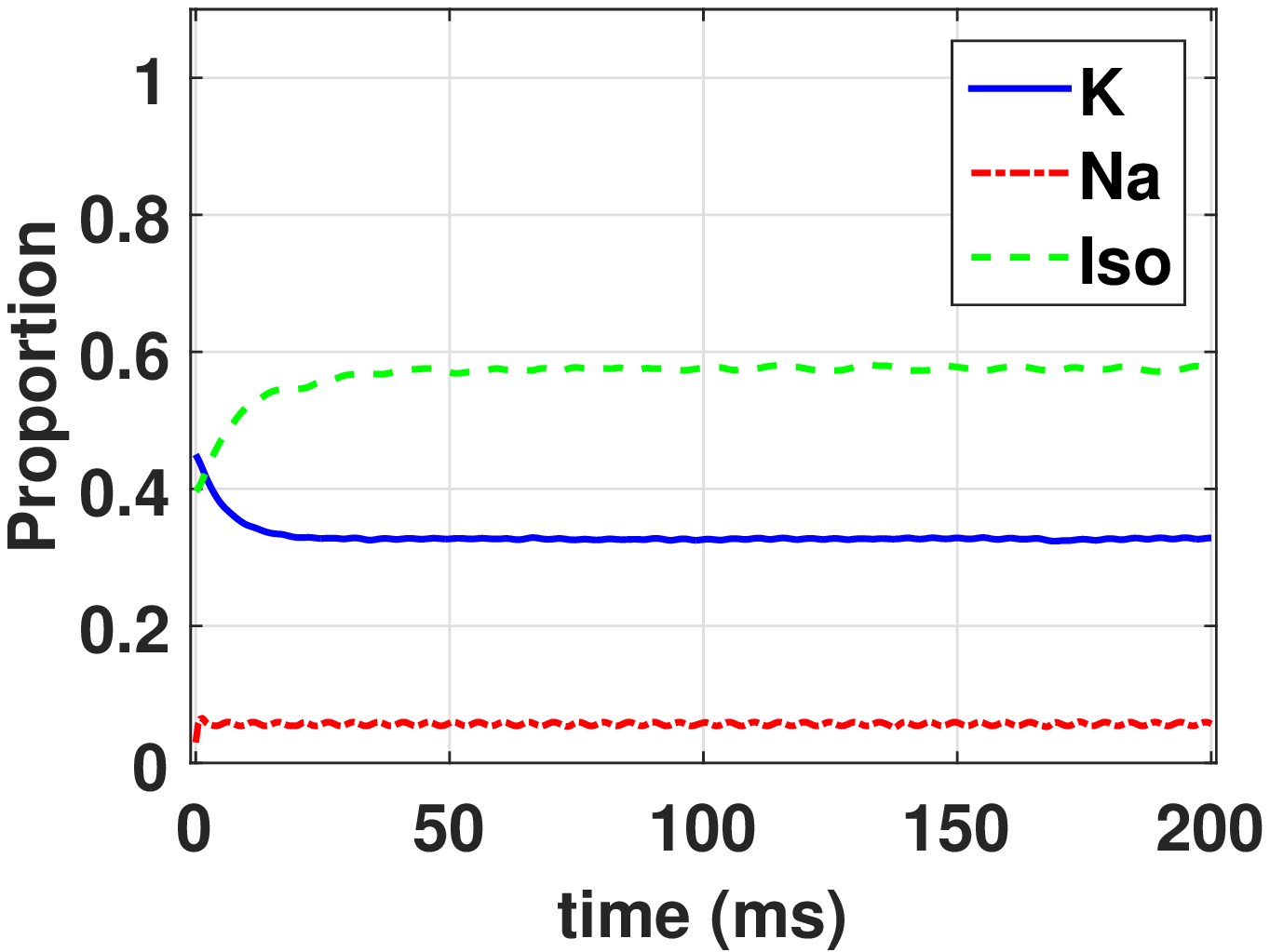} \\
        \hline
    \end{tabular}\\

    {\bf{Fig. 3.}} Graphics $v(t,0)$ (left), $v(t,x)$ (center), $n$, $m$, $h$ (right) for $v_0 = 4.82$, $v^* = 0.5*sin(4/\pi *t)+0.6$, $\rho = 20$
    \end{center}

  	\begin{center}
    \begin{tabular}{|c|c|c|}
    \hline
    \includegraphics[width=2.0in]{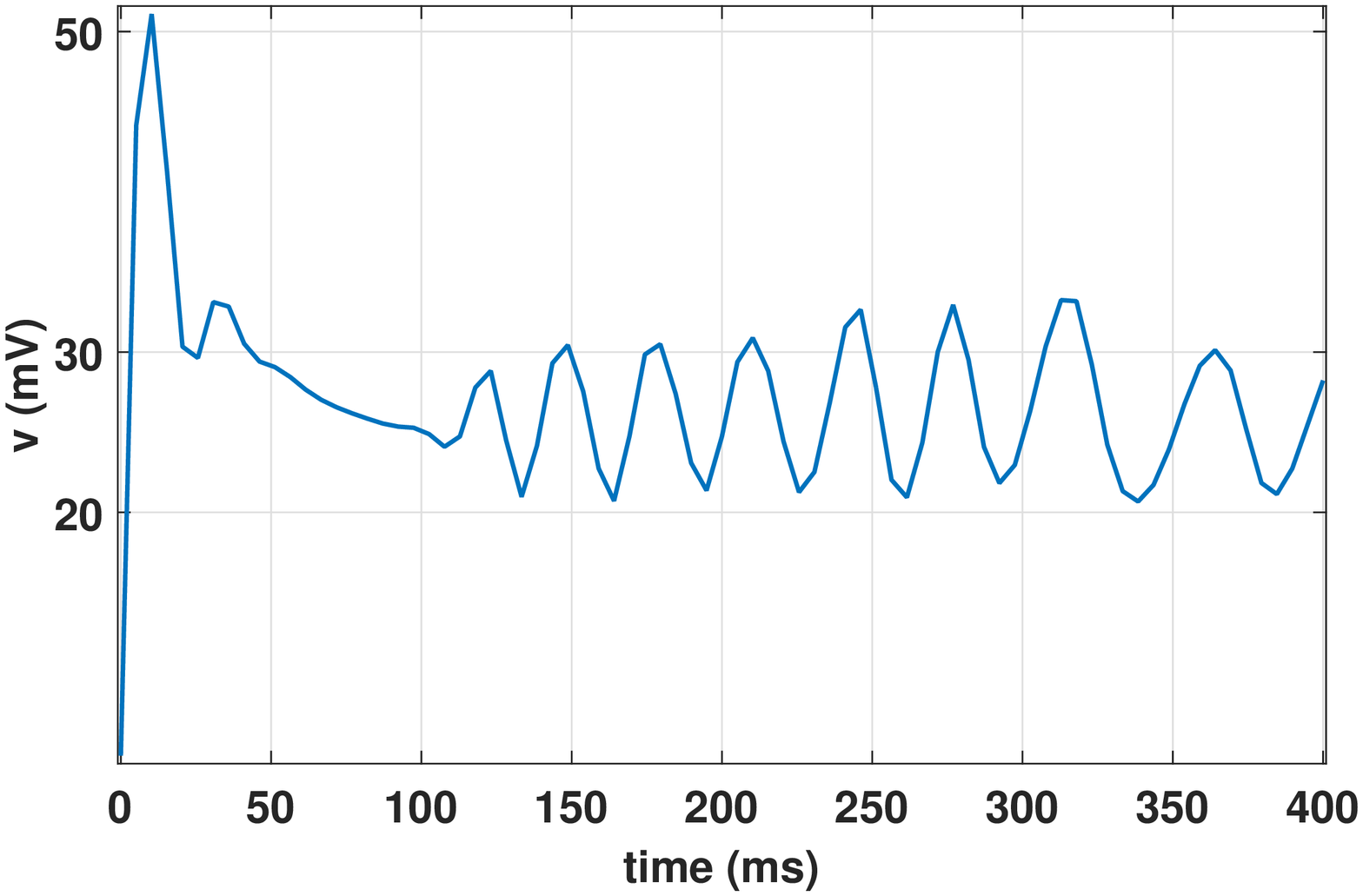} & \includegraphics[width=2.0in]{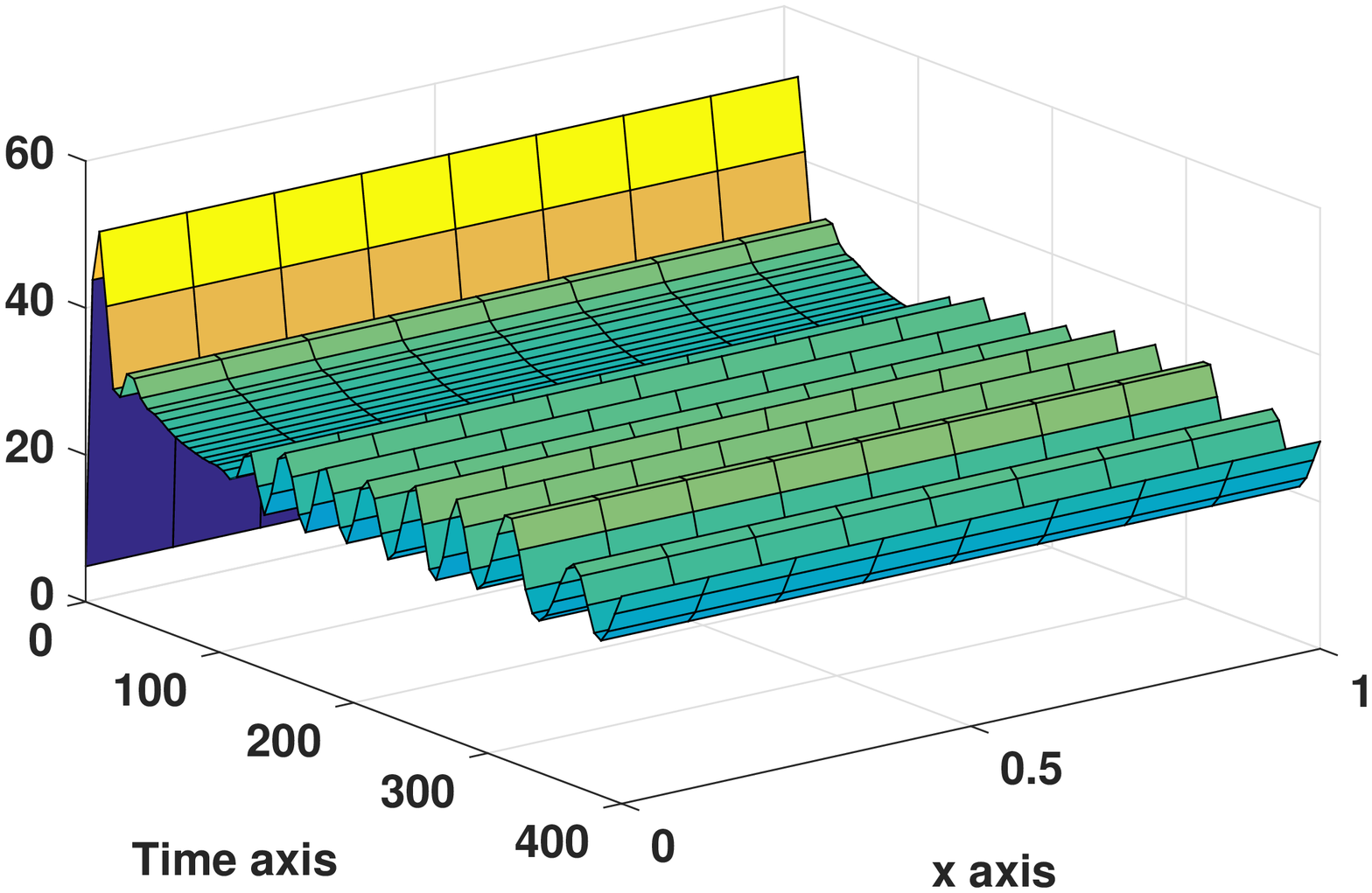} & \includegraphics[width=2.0in]		{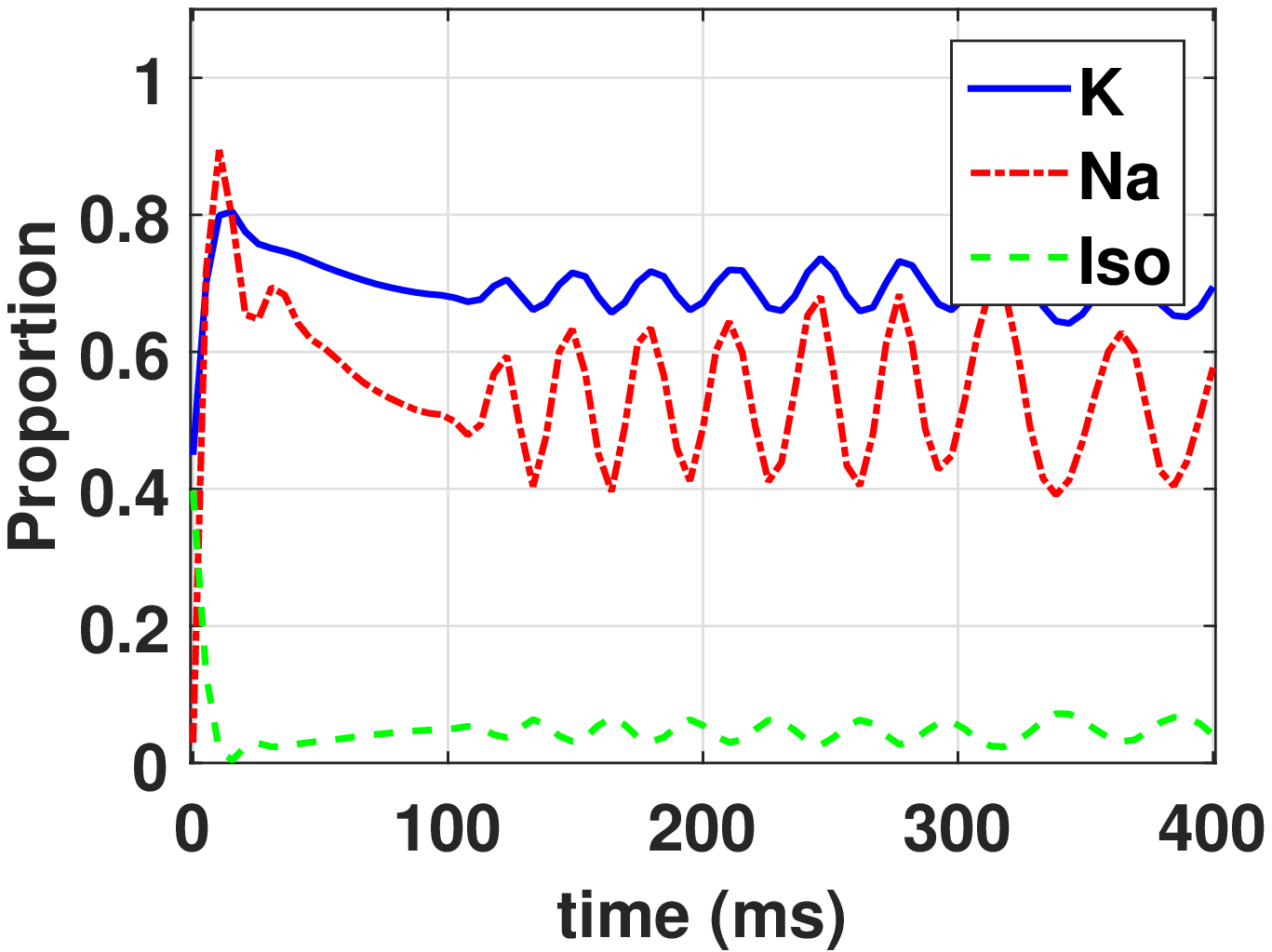} \\
        \hline
    \end{tabular}\\

    {\bf{Fig. 4.}} Graphics $v(t,0)$ (left), $v(t,x)$ (center),$n$, $m$, $h$ (right) for $v_0 = 4.82$, $v^* = 0$, $\rho = 0$, $g_K = 3.8229$
    \end{center}

	 \begin{center}
    \begin{tabular}{|c|c|c|}
    \hline
    \includegraphics[width=2.0in]{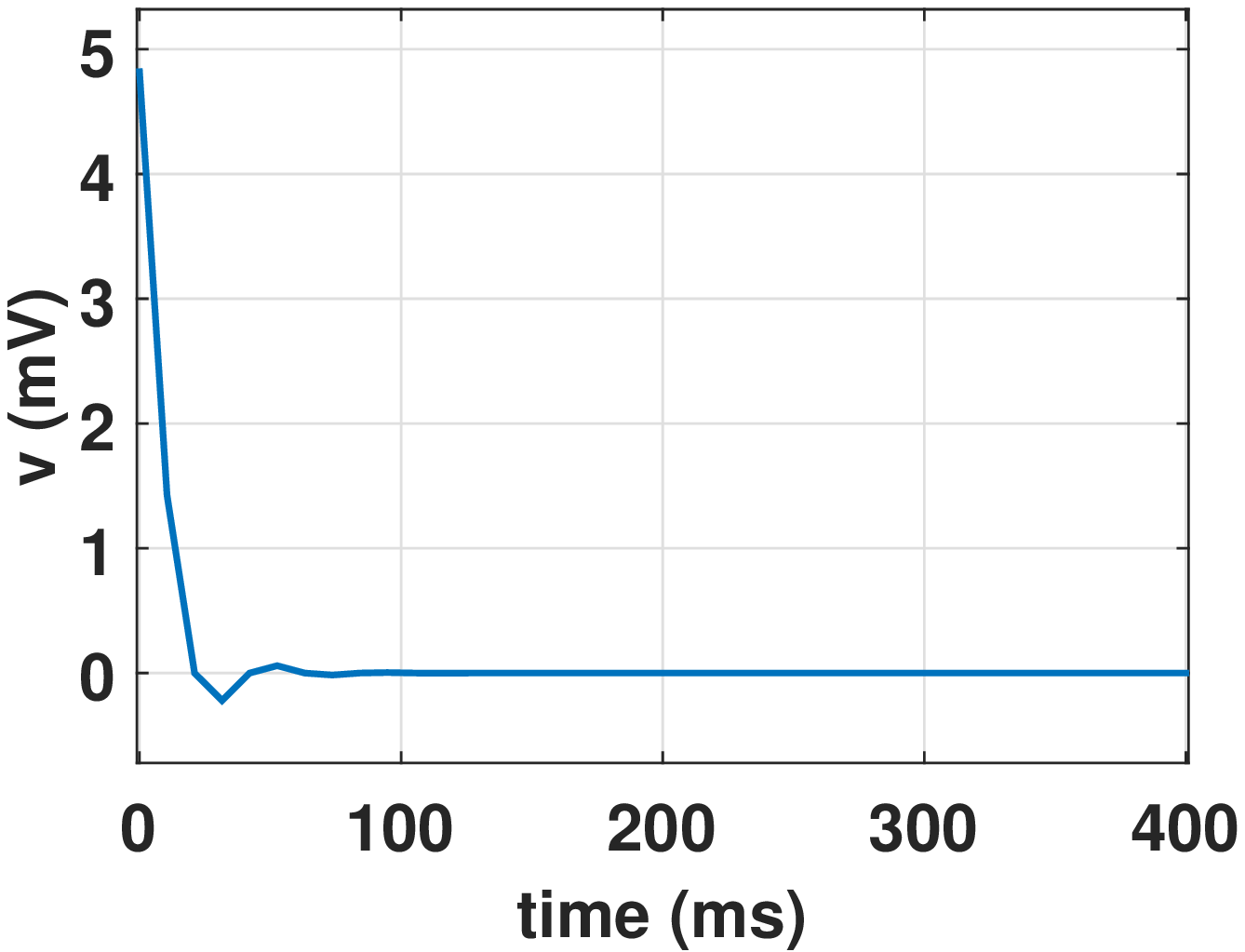} & \includegraphics[width=2.0in]{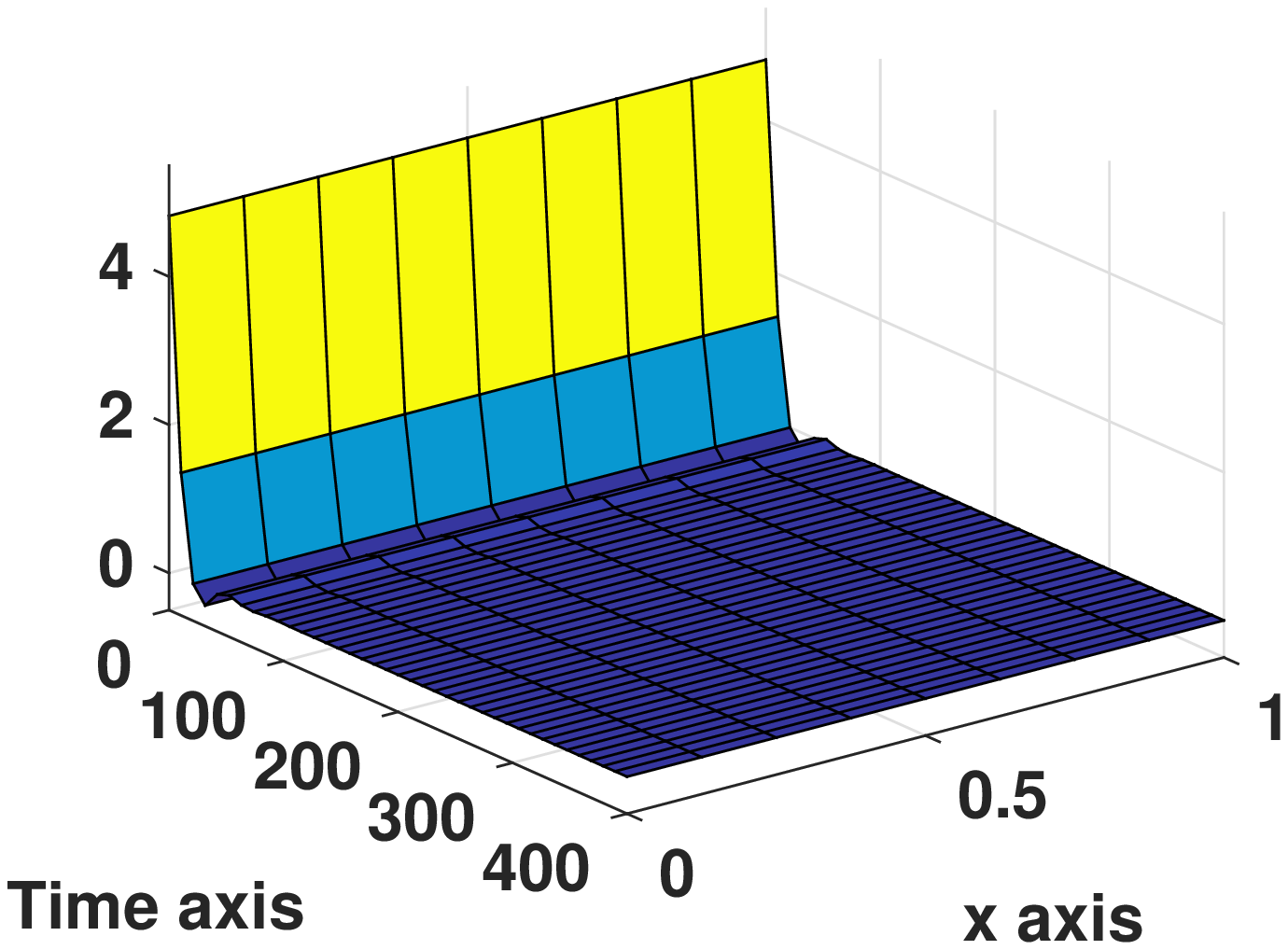} & \includegraphics[width=2.0in]		{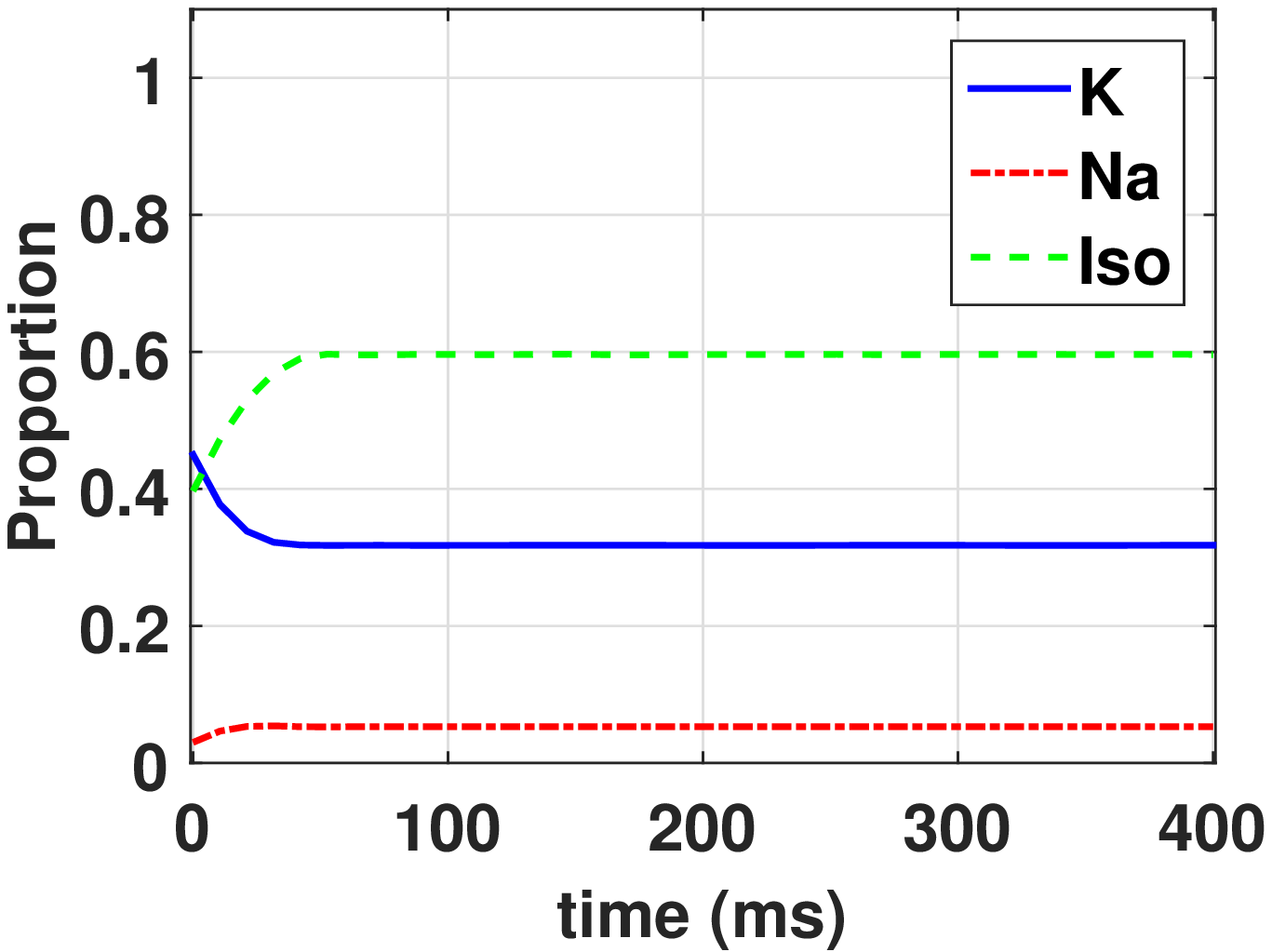} \\
        \hline
    \end{tabular}\\

    {\bf{Fig. 5.}} Graphics $v(t,0)$ (left), $v(t,x)$ (center), $n$, $m$, $h$ (right) for $v_0 = 4.82$, $v^* = 0$, $\rho = 20$, $g_K = 3.8229$
    \end{center}

The final Fig. 6 shows the evolution of the system towards the sliding
stabilization when $v_{0}=0.5\sin (4\pi x)+0.6$ and $g_{K}=36.$ In this case
the graphics $v(t,x_{fixed})$ differ when modifying $x_{fixed}\in \lbrack
0,1]$ due to the space dependence of $v_0$.

 \begin{center}
    \begin{tabular}{|c|c|c|}
    \hline
    \includegraphics[width=2.0in]{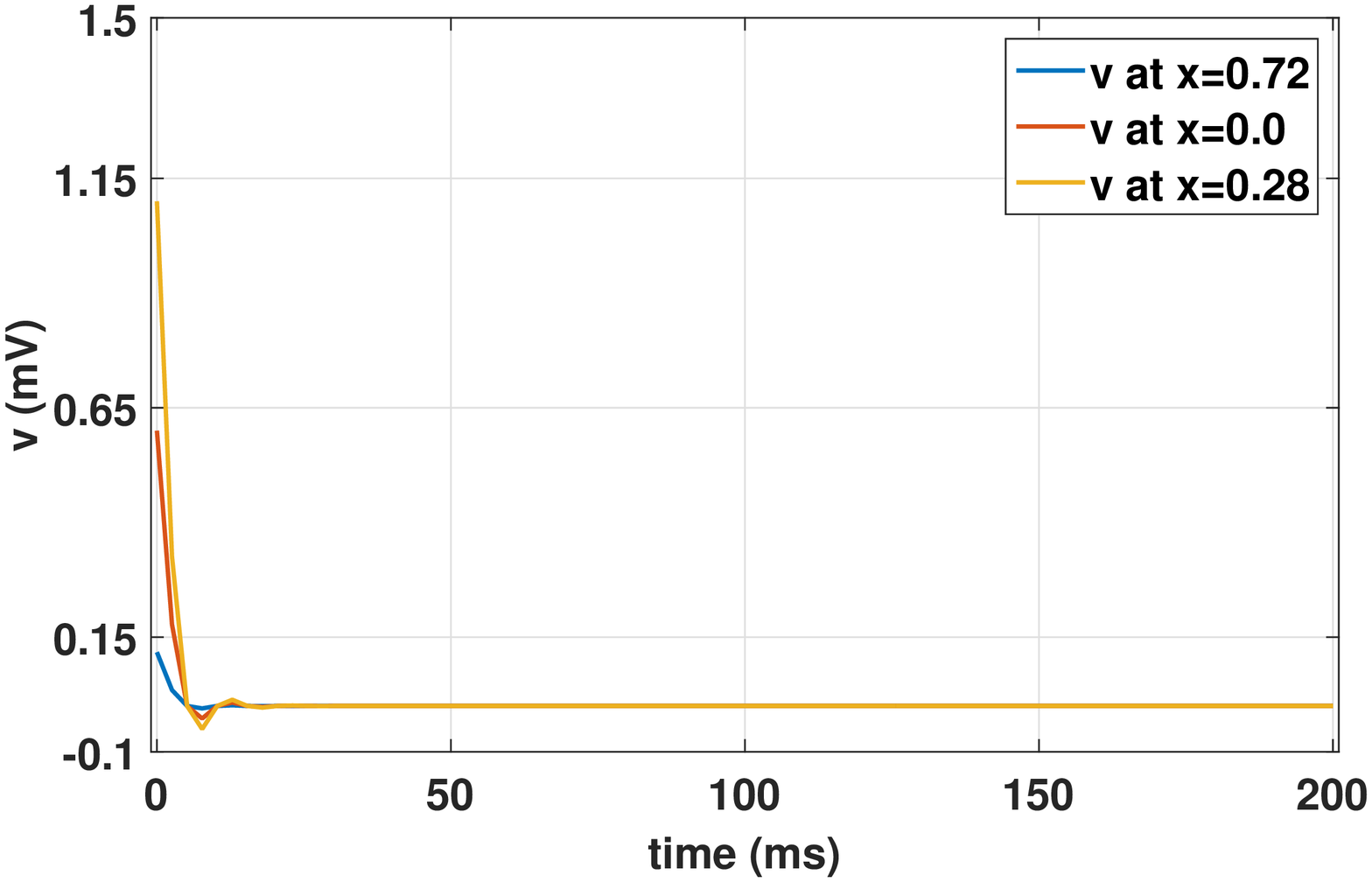} & \includegraphics[width=2.0in]{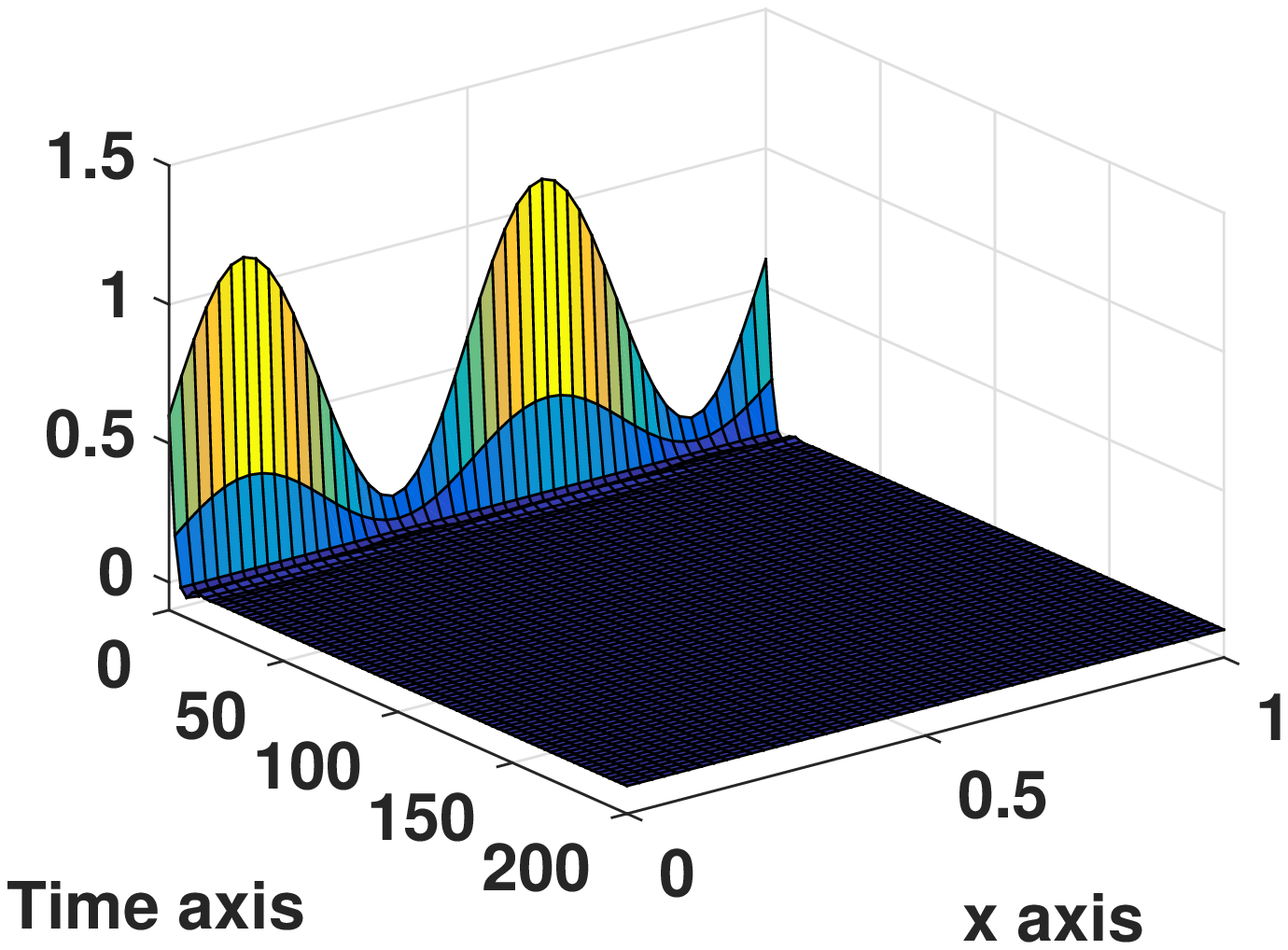} & \includegraphics[width=2.0in]		{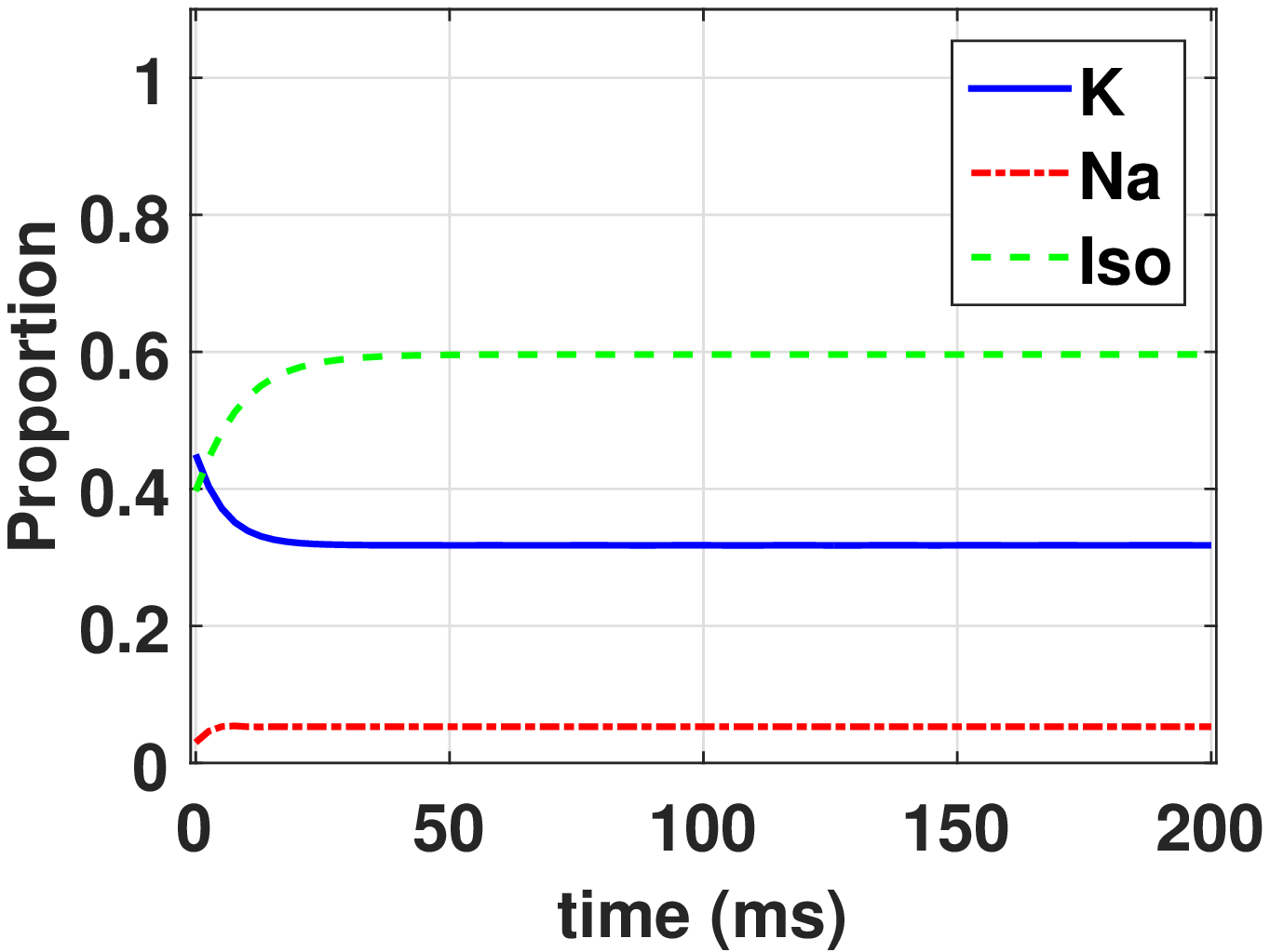} \\
        \hline
    \end{tabular}\\

    {\bf{Fig. 6.}} Graphics $v(t,0)$ (left), $v(t,x)$ (center), $n$, $m$, $h$ at $x=0.5$ (right) for $v_0 = 0.5sin(4\pi x)+0.6$, $v^* = 0$, $\delta =50$, $g_K = 36$, $\rho = 50$
    \end{center}

In all situations starting from a constant $v_{0}$ we observe that the
graphics $v(t,x_{fixed})$ with $x_{fixed}\in \lbrack 0,1]$ are the same,
since the system is invariant to the translation $x\rightarrow x+l$. A
different situation can be observed in Fig. 6 when the initial $v_{0}$
depends on $x.$

While the left and center plots in each figure show the evolution of the
membrane potential, the right ones put into evidence the play between the
other components of the system. The equilibrium potential is determined by
gradients of ionic concentration, through the membrane permeability, and
also, by the effect of the sodium-potassium transport. There is a
concentration of potassium ions inside the cell and a higher concentration
of sodium chloride ions in the external part. At their turn, the
permeabilities of the membrane to sodium and potassium depend on the
membrane potential. The figures on the right show a fast initial inflow of
sodium ions and a subsequent outflow of potassium ions, which define the
action potential generation that follows the stimulation of the
depolarization. The chloride ions do not play their role very well, but they
first exhibit an increase.

We proposed a sliding mode control strategy for the Hodgkin--Huxley model,
by controlling the equation for the membrane potential by a relay type
controller. This permits to reduce the oscillatory movement of the nonlinear
Hodgkin--Huxley system to a stable equilibrium point.

\medskip

\noindent \textbf{Acknowledgments.} This research activity has been
performed in the framework of the Italian-Romanian project \textquotedblleft
Control and stabilization problems for phase field and biological
systems\textquotedblright\ of the Italian CNR and the Romanian Academy,
2017-2019 and was partially supported by a grant of Ministry of Research and
Innovation, CNCS UEFISCDI, Project Number PN-III-P4-ID-PCE-2016-0372, within
PNCDI III, for G.Marinoschi. C. Cavaterra is a member of the Gruppo
Nazionale per l'Analisi Matematica, la Probabilit\`{a} e le loro
Applicazioni (GNAMPA) of the Istituto Nazionale di Alta Matematica (INdAM).
D. En\u{a}chescu is a member of the PhD School of Computer Science,
University of Bucharest.

\end{document}